\documentclass[11pt,fleqn]{article}
\usepackage{amsmath}
\usepackage{amssymb}
\usepackage[pdftex,colorlinks=true,linkcolor=black,%
anchorcolor=black,citecolor=black,urlcolor=black,%
pdftitle={stochvol15June2007},pdfauthor={B. Hanzon and W.
Scherrer}]{hyperref}
\usepackage{graphicx}

\pdfoutput=1

\newtheorem{lemma}{Lemma}[section]

\newtheorem{proposition}{Proposition}[section]

\newcommand{\proof}{{\noindent \bf Proof:\ }}
\newcommand{\eproof}{\hfill $\Box$}

\newcommand{\Rset}{\mathbb{R}}

\newcommand{\Cset}{\mathbb{C}}

\newcommand{\Nset}{\mathbb{N}}

\newcommand{\Ex}{\mathbb {E}}

\newcommand{\ssm}[4]{\left[ \begin{array}{c|c}
{#1 }  &  {#2 }\\
\hline {#3 }  &  {#4 }
\end{array}\right] }

\newcommand{\sssm}[9]{\left[ \begin{array}{cc|c}
{#1 }  &  {#2 }  &  {#3 }\\
{#4 }  &  {#5 }  &  {#6 }\\
\hline {#7 }  &  {#8 }  &  {#9 }
\end{array}\right] }

\newcommand{\rank}{\mathop{\mathrm{rank}}}

\makeatletter
\def\adots{\mathinner{\mkern2mu\raise\p@\hbox{.}
\mkern2mu\raise4\p@\hbox{.}\mkern1mu
\raise7\p@\vbox{\kern7\p@\hbox{.}}\mkern1mu}}
\def\blfootnote{\xdef\@thefnmark{}\@footnotetext}
\makeatother

\begin{document}

\title{Filtering and estimation in stochastic volatility models
 with rationally distributed disturbances}

\author{
Bernard Hanzon \\
School of Mathematical Sciences\\
University College Cork, Ireland\\
\href{mailto:b.hanzon@ucc.ie}{b.hanzon@ucc.ie}\\
\and\\
Wolfgang Scherrer\\
Institut f\"{u}r Wirtschaftsmathematik\\
Technische Universit\"{a}t Wien, Austria\\
Argentinierstr. 8/E105\\
A 1040 Vienna, Austria\\
tel: +43.1.58801-11946\\
\href{mailto:Wolfgang.Scherrer@tuwien.ac.at}{Wolfgang.Scherrer@tuwien.ac.at}
}
\maketitle
%\newpage

\begin{abstract}This paper deals with the filtering problem
for a class  of discrete time stochastic volatility models in
which the disturbances have rational probability density
functions. This includes the Cauchy distributions and Student
t-distributions with odd number of degrees of freedom. Using state
space realizations to represent the rational probability density
functions we are able to solve the filtering problem exactly.
However the size of the involved state space matrices grows
exponentially with each time step of the filter. Therefore we use
stochastically balanced truncation techniques to approximate  the
high order rational functions involved. In a simulation study we
show the applicability of this approach. In addition a simple
method of moments estimator is derived.
\end{abstract}
\begin{description}
\item[Keywords:] stochastic volatility, filtering, rational probability
density function, state space realization, stochastically balanced truncation.
\end{description}
\newpage
%%%%%%%%%%%%%%%%%%%%%%%%%%%%%%%%%%%%%%%%%%%%%%%%%%%%%%%%%%%%%%%%%%%%%%%%%%%%
\section*{Introduction\label{sec:introduction}}
%%%%%%%%%%%%%%%%%%%%%%%%%%%%%%%%%%%%%%%%%%%%%%%%%%%%%%%%%%%%%%%%%%%%%%%%%%%%

In the area of financial time series the Black-Scholes model is
often used for modelling the behaviour of the price of stocks,
exchange rates and other financial time series. This is also the
basis for much of the literature on pricing of derivative
financial instruments such as options. However it is considered to
be a well-known fact that although the volatility is assumed to be
constant in the Black-Scholes model,  in practice it is varying.
This has led to the investigation of more general models in which
the volatility is allowed to vary. One can broadly distinguish
between two types of generalizations. One is the type of model in
which the volatility is varying over time and its dynamic
behaviour is described by some stochastic process. A problem with
such models is that it is generally difficult to solve the {\em
volatility estimation problem} for such models: the calculation of
the conditional density of the volatility at some point in time,
given the observations up till that same point in time, is usually
a difficult task for which there are no closed form expressions.
In the literature there are several proposals to approximate the
conditional density, cf. e.g. \cite{HarveyRuizShephard94},
\cite{MahieuSchotman98}, \cite{BrigoHanzon98}. The other, second
type of model that is used is the ARCH model and its
generalizations (\cite{Engle82}, see also e.g.
\cite{Gourieroux97}), as applied to financial time series. These
models have the advantage that the volatility is again time
varying, and the conditional volatility (also called conditional
heteroskedasticity in this context) is in fact prescribed by the
model as a deterministic function of the past observations. By
construction the problem of estimating the stochastic volatility
has been solved in these models. However one could argue that this
is at the expense of a less transparent model for the underlying
data generating process. In the present paper a model of the first
type will be presented, however with the advantage that for this
model the volatility estimation problem can be solved, as we will
show. Apart from the volatility to be time-varying another feature
of financial time series that is often reported is that it has
{\em fat tails}. In the literature there are many studies that try
to deal with this phenomenon by specifying non-Gaussian
disturbances. This goes back to the work of~\cite{Mandelbrot63}
who suggested to consider the class of stable distributions as
possible distributions for the disturbances. An important example
of stable distributions is given by the Cauchy distributions. More
recent studies seem to favor other distributions, including
Student t-distributions (cf e.g. \cite{CampbellLoMacKinlay97} p.
19, \cite{Lucas96}). In the approach followed in the present paper
all disturbances are allowed which have a {\em rational
probability density function} on the real line. This includes the
Cauchy distributions and Student t-distributions with odd number
of degrees of freedom. In fact it is well-known that the Gaussian
distribution can be approximated by a Student t-distribution of
sufficiently high number of degrees of freedom. Therefore in a
sense the corresponding Gaussian model is a limiting case of the
class of models presented here. It should perhaps be stressed from
the start that there is a price to be paid in the form of high
complexity if one wants to use rational densities of higher
(McMillan) degree. From the point of view of complexity in fact
the estimation problem is easiest when the disturbances have
Cauchy density. In a previous paper a matrix calculus was
developed for performing various calculations with rational
probability density functions (\cite{HanzonOber01}) and applied to
a filtering problem for a class of linear dynamical models. Here
we extend this calculus and show how it can be fruitfully applied
to the non-linear filtering problem of volatility estimation, in a
specific class of stochastic volatility models.

The main extension of the calculus concerns a state-space formula
for the {\em composition} of a proper rational function (which can
be allowed to be a proper rational matrix function) with a proper
rational function, under some minor condition that is required to
ensure the resulting (rational) function is again proper.
%This
%result is new as far as the authors are aware. The state-space
%formula for the composition is of interest in its own right.

It is shown that the conditional probability density functions for
the state are all rational functions in this model class and we
provide an explicit way to calculate these and hence solve the
filtering problem exactly. However as the complexity of the
resulting rational probability density functions increases very
quickly over time, the exact filter cannot be implemented
practically, except during a short period of time. An important
innovation in this respect is the application of an approximation
method stemming from stochastic systems theory, called the SBT
(stochastically balanced truncation) method. This method allows to
find a lower order positive rational density function which
differs at each point on the real line by at most a given
prescribed percentage of the original rational density function.
In an application we use a tolerance level of $2 \%$ giving
excellent results. (The bound used is well-known in stochastic
systems theory and is based on the deep and elegant theory of
Hankel norm approximation). In the implementation of the filter
one needs to switch between various representations of the
rational probability density functions. Numerically reliable
methods are presented to perform these steps. The possibility to
implement the various theoretical ideas in a numerically stable
way is crucial for the success of the practical implementation and
forms one of the key contributions of this paper to the practical
usage of rational probability density functions in filtering
problems. We provide the results of some applications to simulated
data and to empirical FX (foreign exchange) data and present a
number of conclusions. A number of the technical results used are
collected in an appendix.

%%%%%%%%%%%%%%%%%%%%%%%%%%%%%%%%%%%%%%%%%%%%%%%%%%%%%%%%%%%%%%%%%%%%%%%%%%%%
\section{The model class\label{sec:modelclass}}
%%%%%%%%%%%%%%%%%%%%%%%%%%%%%%%%%%%%%%%%%%%%%%%%%%%%%%%%%%%%%%%%%%%%%%%%%%%%

Stochastic volatility models that we will consider are of the following
form:
\begin{equation}\label{eq:model}
\begin{array}{rcl}
X_{t+1}&=&a X_t+ W_{t} \\
Y_t    &=& V(X_t)U_t
\end{array}
\end{equation}
where for each $t\in \Nset=\{1,2,\ldots\},$ the random variables
$X_t, W_t, Y_t, U_t$ take their values in the real numbers, and
where $V(x)$ is a real-valued, positive polynomial function of  $x
\in \Rset;$ $\{W_t,~t\in \Nset\}$ and $\{U_t,~t \in \Nset\}$ are
sequences of jointly stochastically independent real valued random
disturbances with time-invariant probability density functions:
for each $t \in \Nset,$ $W_t$ has rational probability density
function $p_{W}(w),$ $U_t$ has rational  probability density
function $p_{U}(u).$ The initial state $X_1$ has rational
probability density function $p_{X_1}(x).$  The parameter $a$ is a
real number that will be assumed to be unequal to zero for ease of
exposition. In financial applications, the $Y_t$ usually stand for
the returns $Y_t=\log(S_{t+1}/S_{t})$ of some price process
$\{S_t,~t\in \Nset\}.$

A number of remarks can be made about this model class.

\begin{enumerate}
\item The family of rational probability density functions is a very
rich class.
It contains the stable class of Cauchy densities, it contains the
Student distributions
with odd number of degrees of freedom. Under relatively mild conditions,
probability density functions can be approximated by rational
probability density functions,
as follows from results of the theory of rational approximation.
It is well-known that the
Gaussian probability density functions can be approximated for
example by the Student
t-distribution of sufficiently high number of degrees of freedom,
therefore the
Gaussian case appears in a certain sense as a limiting case
of our model class.

\item The function $V$ is a  positive polynomial, i.e. for all
$x\in\Rset,~V(x)> 0.$ Here this is required for technical reasons.
In the literature one finds  other positive functions as
specifications for $V$ as well, for example an exponential
function $V(x)=\exp(\frac{x+\gamma}{2})$ (cf. e.g.
\cite{Taylor86}). If desired one can approximate the exponential
function on any given finite interval by a positive polynomial.
Generalization of the results presented here, to the case in which
$V$ is a non-negative polynomial, i.e. for all
$x\in\Rset,~V(x)\geq 0$ is straightforward.

\item The parameters in the model as well as in the rational probability
density functions of $W_t$ and $V_t,~t=1,2,\ldots,$ are assumed to be
constants here. However they could be taken time-varying if desired. The
resulting filter equations for that case form a straightforward extension
of the filter equations presented in this paper.
\end{enumerate}

%%%%%%%%%%%%%%%%%%%%%%%%%%%%%%%%%%%%%%%%%%%%%%%%%%%%%%%%%%%%%%%%%%%%%%%%%%%%%%%
\section{The filter\label{sec:filter}}
%%%%%%%%%%%%%%%%%%%%%%%%%%%%%%%%%%%%%%%%%%%%%%%%%%%%%%%%%%%%%%%%%%%%%%%%%%%%%%%

We consider the following nonlinear filtering problem: Estimate at
each time $t \in \Nset$ the volatility $V(X_t)$ from the sequence
of observations $Y_1^t:=\{Y_s,~s\in \Nset, s \leq t\}.$ (Note that
we will use the same symbols $Y_1^t:=\{Y_s,~s\in \Nset, s \leq
t\}$ for the random variables and their observed values. This is
to avoid complicating the notation any further. The interpretation
of the symbols as random variables or observed values should be
clear from the context). The solution of such a problem consists
of finding for each $t \in \Nset$ the conditional probability
density function of $X_t$ given $Y_1^t,$ and deriving the desired
estimate of $V(X_t)$ from this. Let the conditional density of
some random variable $Z$ given $Y_1^t$ be denoted by
$p_{Z|Y_1^t}$.

The filter consists of a set of recursive equations by which one
can calculate the conditional probability density function of the
state $X_t$ given the observations $Y_1^t.$
The filter consists of a prediction step and an update
step. In the prediction step one calculates the
conditional density $p_{X_{t+1}|Y_1^t}$ of $X_{t+1}$ given the
observations $Y_1^t$ starting from the conditional
density $p_{X_t|Y_1^t}$ of $X_t$ given $Y_1^t:$
\[
p_{X_{t+1}|Y_1^t}=p_{aX_t|Y_1^t}\star p_{W}.
\]
Here $\star$ denotes convolution.

In the update step one calculates the conditional density of $X_t$ given the
observations $Y_1^t$ from the observation $Y_t$ and the conditional density
of $X_t$ given $Y_1^{t-1},$ using Bayes' rule. Suppose  the
conditional probability density function $p_{X_t|Y_1^{t-1}}$ of
$X_t$ given $Y_1^{t-1}$ is known and the observation $Y_t$ becomes
available. The joint density of $(X_t, Y_t)$ can be obtained from
the joint density of $(X_t, U_t)$   by a change of variables:
\[
 \left(\begin{array}{c} X_t \\ Y_t \end{array} \right)=
 \left(\begin{array}{c} X_t \\ V(X_t)U_t \end{array} \right).
\]
The inverse Jacobian determinant of this change of variables is
$\frac{1}{V(X_t)},$ which is positive because $V$ is a
positive polynomial.  It follows that the joint density of
$(X_t,Y_t)$ is given by
$ p_{X_t,Y_t|Y_1^{t-1}}(x,y)=p_{X_t,U_t|Y_1^{t-1}}(x,\frac{y}{V(x)})
\frac{1}{V(x)}=
p_{X_t|Y_1^{t-1}}(x) p_{U}(\frac{y}{V(x)}) \frac{1}{V(x)}.$
Substituting $y=Y_t,$ we obtain the following expression for the
density of $X_t|Y_1^{t}:$
\[
      p_{X_t|Y_1^{t}}(x)=\frac{1}{c_t}
   p_{X_t|Y_1^{t-1}}(x) p_{U}(\frac{Y_t}{V(x)}) \frac{1}{V(x)},
\]
where
\[
c_t=\int_{-\infty}^{\infty} p_{X_t|Y_1^{t-1}}(x)
p_{U}(\frac{Y_t}{V(x)}) \frac{1}{V(x)} dx.
\]
Since $p_{X_t,Y_t|Y_1^{t-1}}(x,Y_t) =
p_{X_t|Y_1^t}(x)p_{Y_t|Y_1^{t-1}}(Y_t)$ it follows that
\\$c_t=\int p_{X_t,Y_t|Y_1^{t-1}}(x,Y_t) dx =
p_{Y_t|Y_1^{t-1}}(Y_t)$. Therefore we may evaluate the likelihood
as
\[
p_{Y_1,Y_2,\ldots,Y_T}(Y_1,Y_2,\ldots,Y_T) = c_1c_2\cdots c_T.
\]
Note that the value of the normalization constants $c_1,
c_2,\ldots, c_T$ easily follows from Proposition~\ref{prop:Ex}
without the need for explicit integration.

As shown in \cite{HanzonOber01} the convolution of two rational
density functions is a rational function too. Therefore it follows
easily that the conditional density functions defined above will
all be rational, given our assumptions! A way to implement the
filter using ideas from system theory will be presented in the
next sections.

%%%%%%%%%%%%%%%%%%%%%%%%%%%%%%%%%%%%%%%%%%%%%%%%%%%%%%%%%%%%%%%%%%%%%%%%%%%%%
\section{State-space calculus for rational probability density
functions \label{sec:ratdenoperations}}
%%%%%%%%%%%%%%%%%%%%%%%%%%%%%%%%%%%%%%%%%%%%%%%%%%%%%%%%%%%%%%%%%%%%%%%%%%%%%
\subsection{Introduction to the state-space calculus}

 The key idea is to identify rational densities with
spectra of linear, dynamic, continuous time, finite dimensional
systems. This allows us to use concepts and methods from systems
theory; for an overview cf, e.g, \cite{Rugh96}.

Consider a rational non-normalized probability density function
$\rho(x)$. With it we associate a rational function $\Phi(s)$ on
the complex plane which is specified on the imaginary axis by
\[
\rho(x)=\Phi(ix), \;\; \forall x \in \Rset.
\]
Note that $\Phi(\cdot)$ is a rational function which is
nonnegative and integrable on the imaginary axis. Such a function
will be called an \emph{integrable spectral density} in this
paper. The function $\Phi$ has a representation as
\begin{equation}\label{eq:Phipoly}
\Phi(s)=\frac{g_{0}+g_{1}s+\cdots+g_{2q}s^{2q}}{f_{0}+f_{1}s+\cdots+f_{2n}s^{2n}};
\; n>q;\; g_k, f_k\in \Cset
\end{equation}
with coprime polynomials $g(s)=g_{0}+g_{1}s+\cdots+g_{2q}s^{2q}$
and $f(s)=f_{0}+f_{1}s+\cdots+f_{2n}s^{2n}$. Since $\Phi(s)$ is
strictly proper there exists a state space representation, i.e. a
triple $[\tilde{F}\in\Cset^{2n\times 2n}, \tilde{G}\in
\Cset^{2n\times 1}, \tilde{H}\in \Cset^{1\times 2n}]$, such that
\[
\Phi(s) = \tilde{H}(sI_{2n}-\tilde{F})^{-1}\tilde{G}.
\]
Note that we need complex valued triples as $\rho$ may be
non-symmetric. Further note that this representation is not
unique. A state-space transformation
$[\tilde{F},\tilde{G},\tilde{H}]\mapsto
[T\tilde{F}T^{-1},T\tilde{G},\tilde{H}T^{-1}],$ where $T\in
\Cset^{2n\times 2n}$ is a non-singular matrix, leads to a usually
different state-space representation of the same function $\Phi.$
As a shorthand notation for such a state space realization we will
write:
\begin{equation}
\Phi=\pi\ssm{\tilde{F}}{\tilde{G}}{\tilde{H}}{0},\label{eq:ssPhi}
\end{equation}
where $\pi$ will be used in general to denote the mapping that
maps a partitioned matrix $\ssm{A}{B}{C}{D}$ to the corresponding
rational function $C(sI-A)^{-1}B+D.$ It is assumed that the
partitioning involved will be clear from the context in all cases.
Clearly $\pi$ is invariant under state-space transformation.

Since $\Phi(ix)\geq 0$ holds for all $x\in\Rset$, there exists an additive
as well as a multiplicative decomposition of $\Phi(s)$ of the form
\[
\Phi(s) = Z(s)+Z^{*}(s) = K(s)K^{*}(s)
\]
where $Z$ and $K$ are strictly proper.
The rational transfer function $Z(s)$ is called a \emph{spectral summand}
and $K(s)$ is a \emph{spectral factor}.
Here for a rational complex function $G(s)$,
$G^*(s)$ is defined as $G^*(s)=\overline{G(-\bar{s})}$, where $\bar{z}$ denotes
complex conjugation. In particular note that $\Phi^*(s)=\Phi(s)$ holds.
Since $\Phi(s)$ has no poles on the imaginary axis, a stable summand $Z(s)$,
and a stable factor $K(s)$ may be chosen, i.e. $K(s)$ and $Z(s)$
have no pole in the closed right half plane. From now on we always impose
stability on $Z$ and $K$.

Since $Z$, $K$ are strictly proper rational functions, there exist state space representations:
\[
K=\pi\ssm{A}{B}{C}{0}\;;\; K^{*}=\pi\ssm{-A^{*}}{C^{*}}{-B^{*}}{0}
\]
\[
Z=\pi\ssm{A}{M}{C}{0}\;;\; Z^{*}=\pi\ssm{-A^{*}}{C^{*}}{-M^{*}}{0}
\]
Here and elsewhere in this paper $M^*$ denotes the Hermitean
transpose of a matrix $M.$ It is important to note that the two
above realizations may be chosen to share the $A$ and $C$ matrix.

Given these state space realizations for $Z(s)$ and $K(s)$, we may construct
two alternative state space realizations for $\Phi$:
\begin{equation}
\ssm{F}{G}{H}{0}:=\sssm{A}{0}{M}{0}{-A^{*}}{C^{*}}{C}{-M^{*}}{0},
\Phi=\pi\ssm{F}{G}{H}{0}, \label{eq:ssZZ}
\end{equation}
\begin{equation}
\ssm{\bar{F}}{\bar{G}}{\bar{H}}{0}=\sssm{A}{-BB^{*}}{0}{0}{-A^{*}}{C^{*}}{C}{0}{0},
\Phi=\pi\ssm{\bar{F}}{\bar{G}}{\bar{H}}{0}
\label{eq:ssKK}
\end{equation}
using standard formulas for the state space realizations of the
sum and product of two rational functions, see
Appendix~\ref{app:add+mult}.

The \emph{co-degree} of a proper rational function $G(s)$ is
defined as the multiplicity of the zero of $G$ at infinity. Thus
the co-degree of $\Phi$ is $2n-2q$, see~(\ref{eq:Phipoly}).
Clearly the co-degree of $\Phi$ is twice the co-degree of its
spectral factor $K$ and thus is even. For a more detailed
discussion on the co-degree and the zeros of a rational function
see Appendix~\ref{app:codegree}.

As shown in \cite{HanzonOber01} the following proposition
concerning the normalization constant and the moments of a
rational density holds:

\begin{proposition}\label{prop:Ex}
Let $X$ be a real random variable with non-normalized rational
probability density function $\rho$ with corresponding spectral
summand $Z,$ hence $\Phi(ix):=\rho(x)=Z(ix)+Z^*(-ix),$ and let
$(A,M,C)$ be a stable state-space realization of $Z.$ Then $CM$ is
real and positive and $\frac{\rho}{2\pi CM}$ is the {\em
probability} density function corresponding to $X.$ The moments
$\Ex (X^l)$ of $X$ exist for $l=0,\ldots,k-2,$ where $k$ is the
co-degree of $\Phi$ and the moments are given by $\Ex
(X^l)=(-i)^l\frac{CA^l M}{CM},~l=0,1,2,\ldots,k-2.$
\end{proposition}

In  \cite{HanzonOber01} it was shown, i.a., that the operations of
scaling and convolution of rational density functions can be
translated into linear algebra operations on corresponding
state-space realizations. For ease of reference these results are
collected in the Appendix~\ref{app:ratdensop}.

\subsection{The composition formula}

 A key step in the present paper is the construction of
a realization of the rational density function given by
$p_U(\frac{y}{V(x)})/V(x),$ where $y$ is a fixed non-zero real
number and $V(x)$ is a positive polynomial on the real line, if a
realization of $\Phi(ix)=p_U(x)$ is known. Such a realization is
constructed via the following proposition which gives a
realization formula for the composition $G = G_1\circ g_2$,
$G(s)=G_1(g_2(s))$ of two proper rational complex functions $G_1$
and $g_2$. Here we will apply this result to $G_1(ix)=p_U(yx)x$
and $g_2(ix)=i/V(x)$. This implies that we could allow $V(\cdot)$
to be a rational positive function such that $1/V(x)$ is strictly
proper.
%(Also generalization to the case
%where $V$ is a non-negative rational function such that $1/V(x)$
%is strictly proper is possible.)

The only constraint on the pair $G_1$, $g_2$ will be that the
direct feedthrough $d_2=\lim_{s\rightarrow \infty} g_2(s)$ is not
a pole location of $G_1$, i.e. $G_1(d_2)\neq \infty$, because
otherwise the composition $G_1\circ g_2$ would have a pole at
infinity, or in other words it would not be proper rational
function and therefore would not have a state space representation
of the form that we use here. In fact in the proposition we will
allow $G_1$ even to be a rational matrix function, corresponding
to a multi-input, multi-output system in the system theoretic
interpretation.

\begin{proposition}\label{prop:composition}
Let $G_1,g_2$ be proper rational functions with state space
realizations $(A_1\in\Cset^{n_1\times n_1}, B_1\in\Cset^{n_1\times
m_1}, C_1\in\Cset^{p_1\times n_1},D_1\in\Cset^{p_1\times m_1})$
and $(A_2\in\Cset^{n_2\times n_2},b_2\in\Cset^{n_2\times 1},
c_2\in\Cset^{1\times n_2},d_2\in\Cset)$ respectively,
so $G_1(s)=D_1+C_1(sI-A_1)^{-1}B_1$,
$g_2(s)=d_2+c_2(sI-A_2)^{-1}b_2$. Assume that $d_2$ is not an
eigenvalue of $A_1$. Then the composition $G=G_1\circ g_2$ is
again a proper rational (matrix) function with state space
realization $(A,B,C,D)$ given by the formulas
\begin{equation}
\label{eq:comp}
\begin{array}{rcl}
A &=& I_{n_1}\otimes A_2 + (A_1-d_2I_{n_1})^{-1}\otimes b_2c_2
    \in\Cset^{n_1 n_2 \times n_1 n_2}\\
B &=& -(A_1-d_2I_{n_1})^{-1}B_1\otimes b_2
    \in\Cset^{n_1 n_2\times m_1}\\
C &=& C_1(A_1-d_2I_{n_1})^{-1}\otimes c_2
    \in\Cset^{p_1\times n_1 n_2}\\
D &=& D_1-C_1(A_1-d_2I_{n_1})^{-1}B_1
   \in\Cset^{p_1\times m_1}
\end{array}
\end{equation}
\end{proposition}
Here $\otimes$ denotes the Kronecker product,  see
e.g.~\cite{LancasterTismenetsky85}, ch. 12.

\proof
Use will be made by the following inversion formula for
rational matrices that is well known in system theory. Let
$(\tilde{A},\tilde{B},\tilde{C},\tilde{D}),$ with $\tilde{D}$
invertible, denote the state space realization of a proper
rational function
$\tilde{G}(s)=\tilde{D}+\tilde{C}(sI-\tilde{A})^{-1}\tilde{B}$.
Its inverse is given by $(\tilde{G}(s))^{-1}= \tilde{D}^{-1}-
\tilde{D}^{-1}\tilde{C}(sI-\tilde{A}+\tilde{B}\tilde{D}^{-1}\tilde{C})^{-1}\tilde{B}\tilde{D}^{-1}$.

We need to show that the rational matrix $G(s)$ with state space realization given by~(\ref{eq:comp})
is equal to $G_1(g_2(s))$. In order to do that we calculate $G(s)$ as follows:
\[
\begin{array}{rcl}
G(s) &=& D+C(sI-A)^{-1}B\\
&=& D_1-C_1(A_1-d_2I_{n_1})^{-1}B_1 \\
& & \;\;\;\;    -\left(C_1(A_1-d_2I_{n_1})^{-1}\otimes c_2\right)
    H^{-1}
    \left((A_1-d_2I_{n_1})^{-1}B_1\otimes b_2\right)\\
&=& D_1-C_1(A_1-d_2I_{n_1})^{-1}B_1 \\
& & \;\;\;\;   -C_1(A_1-d_2I_{n_1})^{-1}(I_{n_1}\otimes c_2)
    H^{-1}
    (I_{n_1}\otimes b_2)(A_1-d_2I_{n_1})^{-1}B_1\\
&=& D_1+ C_1\left[-(A_1-d_2I_{n_1})^{-1} \right. \\
& & \;\;\;\; \left. -(A_1-d_2I_{n_1})^{-1}(I_{n_1}\otimes c_2)
    H^{-1}
    (I_{n_1}\otimes b_2)(A_1-d_2I_{n_1})^{-1}\right]B_1.\\
\end{array}
\]
where
\[
\begin{array}{rcl}
H &=& sI_{n_1n_2} - I_{n_1}\otimes A_2 -
        (A_1-d_2I_{n_1})^{-1}\otimes b_2c_2 \\
    &=& sI_{n_1}\otimes I_{n_2} - I_{n_1}\otimes A_2 -
        (I_{n_1}\otimes b_2)(A_1-d_2I_{n_1})^{-1}(I_{n_2}\otimes c_2).
\end{array}
\]
This expression has the form $D_1+C_1(\tilde{G}(s))^{-1}B_1$, where
$\tilde{G}(s)=\tilde{D}+\tilde{C}(sI-\tilde{A})^{-1}\tilde{B}$ and
\[
\tilde{D}=-(A_1-d_2I_{n_1}),\;\tilde{C}=(I_{n_1}\otimes c_2),\;\tilde{A}=I_{n_1}\otimes A_2,\;
\tilde{B}=(I_{n_1}\otimes b_2).
\]
It follows that
\[
\begin{array}{rcl}
\tilde{G}(s)&=&-(A_1-d_2I_{n_1})+(I_{n_1}\otimes c_2)\left(sI_{n_1 n_2}-I_{n_1}\otimes A_2\right)^{-1}
(I_{n_1}\otimes b_2)\\
&=&I_{n_1}\otimes \left(d_2+c_2(sI_{n_2}-A_2)^{-1}b_2\right)-A_1 \\
&=&I_{n_1}\left(d_2+c_2(sI_{n_2}-A_2)^{-1}b_2\right)-A_1.
\end{array}
\]
and thus
\[
G(s)=D_1+C_1\left( I_{n_1}\left(d_2 + c_2 (sI_{n_2}-A_2)^{-1} b_2 \right) -A_1\right)^{-1}B_1
= G_1(g_2(s)).
\]
\eproof

{\em Remark.} A special case of the composition formula can be
found in the theory of phase-type distributions in statistics. See
e.g. \cite{Neuts101} and \cite{FackrellThesis}, equations (2.5.1),
(2.5.2).

\subsection{Transformations between the various representations
of rational density functions\label{sec:rsum2rfac}}

 In order to implement the filter we need various
different representations for the conditional densities, i.e. the
integrable spectral density $\Phi$, the spectral summand $Z$,
$Z+Z^*=\Phi$, and the spectral factor $K$, $KK^*=\Phi$. Thus we
need procedures which compute such a representation from any of
the others. The computation of $\Phi$ from $Z$ or $K$ follows from
the formulas in Appendix~\ref{app:add+mult}. The computation of a
spectral summand $Z$ from given $K$ or $\Phi$ is dealt with in the
appendices~\ref{app:rfac2rsum} and~\ref{app:rpdf2rsum}
respectively. The most demanding task is the computation of a
spectral factor given a spectral summand and this will be
presented in the following subsection. It should be noted that
this spectral factorization problem is a standard problem in
systems theory. However most of the literature deals with the case
where $\lim_{s\rightarrow \infty} \Phi(s)=R>0$ holds, i.e where
there is no zero at infinity. Thus we found we had to develop a
numerically robust procedure for the case where $\Phi(s)$ has a
zero at infinity.

%%%%%%%%%%%%%%%%%%%%%%%%%%%%%%%%%%%%%%%%%%%%%%%%%%%%%%%%%%%%%%%%%%%%%%%%%%%%%
%\subsection{Computation of a spectral factor from a spectral summand\label{sec:rsum2rfac}}
%%%%%%%%%%%%%%%%%%%%%%%%%%%%%%%%%%%%%%%%%%%%%%%%%%%%%%%%%%%%%%%%%%%%%%%%%%%%%

Let a spectral summand $Z(s)=C(sI_{n}-A)^{-1}M$
be given. Now the task is to compute $B$ such that
$K(s)=C(sI_{n}-A)^{-1}B$
is a spectral factor, i.e. such that $\Phi(s)=Z(s)+Z^{*}(s)=K(s)K^{*}(s)$
holds. The basic tool for this conversion is the so called \emph{positive
real lemma}:

\begin{lemma}\label{lem:posreal}
A stable rational function $Z(s)=C(sI_{n}-A)^{-1}M$
is \emph{positive real,} i.e. $Z(ix)+Z^{*}(ix)\geq 0$ for all $x\in\Rset$,
if and only if there exists a solution $P$ \emph{}of the \emph{linear
matrix inequality} (LMI)
\begin{equation}\label{eq:LMI}
L(P)=\left[\begin{array}{cc}
-AP-PA^{*} & M-PC^{*}\\
M^{*}-CP & 0\end{array}\right]\geq0
\end{equation}

If $P$ is a solution, then $\Phi(s)=Z(s)+Z^{*}(s)=K(s)K^{*}(s)$,
where $K(s)=C(sI_{n}-A)^{-1}B$ and $B\in\Cset^{n\times r}$ is
determined from
\begin{equation}\label{eq:rankone}
L(P)=\left[\begin{array}{c}
B\\
0\end{array}\right]\left[\begin{array}{c}
B\\
0\end{array}\right]^{*}.
\end{equation}

\end{lemma}
For a proof of this lemma see e.g.~ \cite{Faurre76}. In addition
we remark:

\begin{enumerate}
\item
    The rank of $L(P)$ determines the column dimension of the
    function $K(s)$. In particular it can be shown that there always
    exist square factors $K$. Since we here deal exclusively with the
    scalar case, we are only interested in \emph{rank one solutions},
    i.e. in solutions $P$ where $\rank L(P)=1$.
\item
    By the asymptotic stability of $A$ it follows that any solution $P$ of the
    LMI is positive semidefinite.
\item The solution set $\mathcal{P}=\{P\,|\, L(P)\geq 0\}$ is convex and bounded.
    If the set $\mathcal{P}$ is non-empty it  contains a minimal and a maximal
    element, $\underline{P}$ and $\overline{P}$ say, i.e.
    $\underline{P}\leq P \leq \overline{P}$ holds for all $P\in \mathcal{P}$.
    The minimum element corresponds to a \emph{minimum phase} factor,
    $\underline{K}$ say, i.e. all zeros of $\underline{K}(s)$ are
    in the closed left half plane:
    $\underline{K}(s)=0$ $\Rightarrow$ $\Re(s)\leq 0$.
    Analogously $\overline{P}$ gives a \emph{maximum phase}
    factor $\overline{K}$,
    i.e $\overline{K}(s)=0$ $\Rightarrow$ $\Re(s)\geq 0$.
\end{enumerate}

By the positive real lemma it follows that the computation of the
spectral factor is equivalent to the solution of the above LMI. A
solution of this LMI will be constructed via the computation of
what is known as a deflating subspace of the $(2n+1)\times (2n+1)$
dimensional pencil:
\begin{equation}\label{eq:pencil}
\lambda E-N :=
\ssm{\lambda I_{2n}-F}{G}{-H}{0}=
\sssm{\lambda I_{n}-A}{0}{M}{0}{\lambda I_{n}+A^{*}}{C^{*}}{-C}{M^{*}}{0}
\end{equation}
Note that the eigenvalues of this pencil are the zeros of
$\Phi(s).$ See Appendix~\ref{app:codegree} for background material
on pencils of the above form!

Suppose for the moment that we have given a (rank one) solution
$P=P^{*}$ of the LMI and the corresponding factor
$K(s)=C(sI_n-A)^{-1}B$. Furthermore let $k=2c$ be the co-degree of
$\Phi$ and thus $c$ is the co-degree of $K$. By some easy algebra
it follows that
\[
\sssm{\lambda I_n-A}{0}{-M}{0}{\lambda I_n+A^*}{-C^*}{-C}{M^*}{0}
\left[\begin{array}{cc} P & 0 \\ I_n & 0 \\ 0 & 1 \end{array}\right] =
\left[\begin{array}{cc} P & B \\ I_n & 0 \\ 0 & 0 \end{array}\right]
\ssm{\lambda I_n+A^*}{-C^*}{B^*}{0}
\]
This implies
\begin{enumerate}
\item
\[
\left[\begin{array}{cc} P & 0 \\ I_n & 0 \\ 0 & 1 \end{array}\right]
\]
is a basis for a deflating subspace of the pencil $(\lambda E-N)$.
\item
\[
\ssm{\lambda I_n+A^*}{-C^*}{B^*}{0}
\]
is the pencil corresponding to the zeros of $K^*(s)$ and thus has
a $(c+1)$ dimensional infinite elementary divisor and an $(n-c)$
dimensional finite divisor corresponding to the finite zeros of
$K^*(s)$.
\end{enumerate}

In order to construct a rank one solution of the
LMI~(\ref{eq:LMI}) we therefore have to compute an
$(n+1)$-dimensional divisor of the pencil~(\ref{eq:pencil}) which
itself has a $(c+1)$ dimensional infinite elementary divisor and
an $(n-c)$ dimensional finite divisor. Let $Z\ \in
\Cset^{(2n+1)\times (n+1)}$ be a basis for the corresponding
deflating subspace, where in addition it is assumed that the first
$c+1$ columns form a basis for the $(c+1)$ dimensional deflating
subspace corresponding to the $(c+1)$ dimensional infinite
elementary divisor. By the discussions in
Appendix~\ref{app:codegree} it follows that $Z$ may be partitioned
as
\begin{equation}\label{eq:Zpart}
Z = \left[\begin{array}{ccc}
0 & Z_{12} & Z_{13} \\
0 & Z_{22} & Z_{23} \\
z_{31} & 0 & z_{33} \\
\end{array}\right],\,
Z_{12}, Z_{22} \in \Cset^{n\times c}, \, Z_{13}, Z_{23} \in
\Cset^{n\times (n-c)}.
\end{equation}
Note that $Z$ can be written as
\[Z=\left[\begin{array}{cc}P&0\\I_n&0\\0&1\end{array}\right]T,\]
where $T$ is a $(n+1)\times (n+1)$ non-singular matrix.
Hence the solution $P$ is obtained from
\begin{equation}\label{eq:ZtoP}
P = [Z_{12},Z_{13}][Z_{22},Z_{23}]^{-1}.
\end{equation}
The only remaining choice is the choice of the finite eigenvalues,
which determine the zeros of the factor $K^*$. E.g. in order to
get the minimal solution $\underline{P}$ (the minimum phase factor
$\underline{K}$) one has to choose the $(n-c)$ anti stable
eigenvalues $\Re(\lambda_i)>0$. On the other hand choosing the
stable eigenvalues $\Re(\lambda_i)<0$ gives the maximum element
$\overline{P}$ and the maximum phase factor $\overline{K}$.

This procedure works provided that there are no zeros on the
imaginary axis (except for the zero at infinity). Therefore for
our implementation of the filter in addition we assume that
$p_U(x)$ and $p_{X_1}(x)$ are strictly positive, which implies
that all conditional densities in the filter will be strictly
positive. However the numerical implementation still may run into
trouble if there are zeros ``close'' to the imaginary axis!

The actual procedure is now as follows: Start with the
pencil~(\ref{eq:pencil}) and bring it to the staircase
form~(\ref{eq:staircase2},\ref{eq:staircase3}). Apply a QZ
transformation to the lower right $((2n-2c)\times (2n-2c))$
dimensional block to bring the whole pencil into a QZ form, see
e.g. \cite{GolubVanLoan89}. So $QEZ$ and $QNZ$ are upper
triangular matrices and $Q$ and $Z$ are both unitary. Next by a
sequence of $2\times 2$ orthogonal transformations the diagonal
elements corresponding to the $n-c$ anti-stable (stable)
eigenvalues are shifted to positions $c+2,c+3,\ldots,n+1$, without
losing the triangular structure. The desired basis for the
deflating subspace then is given by the first $n+1$ columns of the
final $Z$ matrix. Finally compute $P$ as described
in~(\ref{eq:Zpart}, \ref{eq:ZtoP}) and $B$ from ~(\ref{eq:LMI}),
(\ref{eq:rankone}) in Lemma \ref{lem:codegG}.

\subsection{Description of the filter in terms of state space
formulas}

 We can now describe how the filter could be calculated
using state space formulas. Recall that for each rational
probability density we associate three rational functions, namely
the spectral density $\Phi$, the spectral summand $Z$ and the
spectral factor $K$. Above it has been discussed how one can
obtain the state space realization of $K$ given the realization of
$Z$. In Appendix~\ref{app:rfac2rsum} it is shown how to compute a
state space realization of $Z$ given a state space realization of
$K$ and in~\ref{app:rpdf2rsum} a realization of $Z$ is computed
from a realization of $\Phi$. A state space realization of $\Phi$
given a state space realization for $Z$ or $K$ follows from the
formulas in Appendix~\ref{app:add+mult}. Therefore we can switch
between these state space realizations as needed.

To start consider the probability density $p_{X_t|Y_1^{t-1}}$ or
for $t=1$ the density $p_{X_1}$. Calculate its spectral factor,
$K_1$ say. Consider the spectral density $\Phi_U$ of $U$ and
construct a state space realization for $G(s)=\Phi_U(-i Y_t s)
(-is)$, where $Y_t$ is the observed output variable. Construct a
state space realization of the spectral density function
$g(s)=i/V(-is)$. Then use the composition formula to obtain a
realization of the spectral density $G\circ g$ of
$p_U(Y_t/V(x))/V(x)$. Calculate the realization of the spectral
factor, $K_2$ say, of this density. Construct the product of $K_1$
and $K_2$ (see Appendix~\ref{app:add+mult}), this gives the
realization of the spectral factor of $c_t p_{X_t|Y_1^t}$.
Calculate the corresponding spectral summand. Compute the
normalization factor $c_t$ from Proposition~\ref{prop:Ex} and
compute the realization of the spectral summand of
$p_{X_t|Y_1^t}$. This will be the input for the prediction step.

Calculate the realization of the spectral summand of
$p_{aX_t|Y_1^t}$ by using the formula for the scaling, see
Proposition~\ref{prop:scaling+convolution}, part (i) in Appendix
\ref{app:ratdensop}. Form the realization of the spectral summand
of $p_W$. Construct the realization of the spectral summand of
$p_{X_{t+1}|Y_1^t}=p_{aX_t|Y_1^t}\star p_W$ using the convolution
formula given in Proposition \ref{prop:scaling+convolution}, part
(ii), in Appendix \ref{app:ratdensop}.

Now, as soon as a new observation $Y_{t+1}$ becomes available one
can proceed to a new update step.

The behaviour of the co-degree and the state-space dimension and
McMillan degree of the conditional densities in the filter can now
be described. First consider the co-degree. Let $k_{1|0}$ denote
the co-degree of $p_{X_1}$ and $k_{t|s}$ the co-degree of
$p_{X_t|Y_1^{s}}.$ We know that the co-degree of a rational
probability density is even. Let $d$ denote the degree of the
polynomial $V,$ then the rational function $g$ constructed above
has co-degree $d.$ Because by assumption the probability density
of $U$ is (strictly) positive, it follows that the co-degree of
$G\circ g$ is $d.$ Hence the co-degree of $p_{X_t|Y_1^t}$ is
$k_{t|t}=k_{t|t-1}+d.$ As non-zero scaling does not affect the
co-degree and convolution of two rational densities leads to a
rational density with co-degree equal to the minimum of the two
co-degrees of the arguments of the convolution (see the notes
after Proposition~\ref{prop:scaling+convolution} in Appendix
\ref{app:ratdensop}) we find $k_{t+1|t}=\min(k_{t|t}, k_W),$ where
$k_W$ denotes the co-degree of the rational density of $W.$ This
makes that the co-degrees of the conditional densities
$p_{X_{t+1}|Y_1^t}$, $t\geq 1$ are bounded by $k_W.$ This will
turn out to be important in the next section.

Now let us turn to an analysis of the state-space dimensions and
the McMillan degree. For the composition and convolution of two
rational functions the state dimension of the output is the
product of the respective dimensions of the inputs (see
Propositions~\ref{prop:composition},
\ref{prop:scaling+convolution}) whereas for the product the
dimensions add up, see Appendix~\ref{app:add+mult}. Therefore in
each step the state dimension of the realization of the
conditional densities tends to increase dramatically. (Note that
together with the result on co-degrees this suggests that the
resulting conditional probability density functions will have a
non-trivial numerator, even if one uses Student-t or Cauchy
densities for the disturbances). Theoretically it is possible that
the resulting state-space realization is not-minimal, in which
case the McMillan degree would be smaller than the state-space
dimension and a state-space reduction procedure could be applied.
In practice we do not expect this to happen very often. However if
this is the case {\em approximately}, one can profitably apply
model reduction techniques, to keep the state-space dimensions
manageable. We suggest to apply model reduction at each time step
to approximate the high degree rational density
$p_{X_{t+1}|Y_1^t}$ by a lower degree rational density. This will
be the topic of the next section.

%%%%%%%%%%%%%%%%%%%%%%%%%%%%%%%%%%%%%%%%%%%%%%%%%%%%%%%%%%%%%%%%%%%%%%%%%%%%%%
\section{Balancing and balanced model reduction\label{sec:bal}}
%%%%%%%%%%%%%%%%%%%%%%%%%%%%%%%%%%%%%%%%%%%%%%%%%%%%%%%%%%%%%%%%%%%%%%%%%%%%%%

There are many possibilities for model reduction. The challenge here is that the
 approximant has to be nonnegative on the real axis. In terms of the
corresponding spectral density this means that the approximant
spectral density has to be nonnegative on the imaginary axis. In
terms of the spectral summand this means that the approximant has
to be positive real.

One well known method to achieve this is the so called ``positive
real balanced truncation'' technique, see e.g.~\cite{DesaiPal84},
which we will shortly explain here.

Note that the solution set $\mathcal{P}$ contains two particular elements,
namely the minimal and the maximal element,
$\underline{P} \leq \overline{P}$ say,
and it has been discussed in section~\ref{sec:rsum2rfac}  how to compute
these elements. A state space realization $[A,M,C]$ of a spectral summand
$Z(s)=C(sI_{n}-A)^{-1}M$ is called \emph{positive real balanced} iff
\[
\underline{P}=\overline{P}^{-1}=\Sigma=
\mathrm{diag}(\sigma_{1},\ldots,\sigma_{n}),\;\;
\sigma_1\geq \sigma_{2}\geq \cdots \geq \sigma_n
\]
holds. The $\sigma_{i}$'s are called the \emph{positive real
singular values} of $Z$. Since $\underline{P}\leq \overline{P}$
holds and since the squared singular values $\sigma_i^2$ are the
eigenvalues of $\underline{P}\overline{P}^{-1}$ it follows that
these singular values are bounded by $0 \leq \sigma_i\leq 1$.
Furthermore it is known (see \cite{Green88LAA}, Theorem 4.1) that
$\sigma_1=\cdots=\sigma_c=1$ and $1> \sigma_{j}$ for all $j=c+1,
\ldots, n,$ holds, where $k=2c$ is the co-degree of $\Phi=Z+Z^*$.

It is easy to see that a state space transformation
\newline
$[A,M,C] \rightarrow [TAT^{-1},TM,CT^{-1}]$, where $T\in
\Cset^{n\times n}$ is a non singular matrix, transforms
$\underline{P}$ and $\overline{P}^{-1}$ as
$\underline{P}\rightarrow T\underline{P}T^*$ and
$\overline{P}^{-1}\rightarrow T^{-*}\overline{P}^{-1}T^{-1}$.
Therefore such a balanced realization may be obtained by the
following procedure. Suppose $\underline{P}$ and $\overline{P}$
are given and let
$\underline{P}=\underline{P}^{1/2}\underline{P}^{*/2}$ and
$\overline{P}=\overline{P}^{1/2}\overline{P}^{*/2}$ be some
arbitrary factorization of these positive definite matrices. Here
$M^{1/2}$ denotes a square root of a positive definite matrix
$M\geq 0$, i.e. $M=M^{1/2} (M^{1/2})^*$. In addition we use the
notations $M^{*/2}=(M^{1/2})^*$, $M^{-1/2}=(M^{1/2})^{-1}$ and
$M^{-*/2}=(M^{*/2})^{-1}$. Next let
$\underline{P}^{*/2}\overline{P}^{-*/2}=U\Sigma V^{*}$ with $U,~V$
unitary matrices, be a singular value decomposition (SVD). The
state space transformation
\[
T=\Sigma^{-1/2}V^{*}\overline{P}^{-1/2}=\Sigma^{1/2}U^*\underline{P}^{-1/2}
\]
then gives the desired balanced realization, since
\[
\begin{array}{rcl}
T\underline{P}T^{*} & = &
\Sigma^{1/2}U^*\underline{P}^{-1/2}
\,\underline{P}\,
\underline{P}^{-*/2}U\Sigma^{1/2} = \Sigma \\
T^{-*}\overline{P}^{-1}T^{-1} & = &
\Sigma^{1/2}V^{*}\overline{P}^{*/2}
\overline{P}^{-1}
\overline{P}^{1/2}V\Sigma^{1/2} =
\Sigma
\end{array}
\]
Let $[\bar{A},\bar{M},\bar{C}]$ denote the balanced realization
obtained by this procedure and let these matrices by partitioned
as
\[
\left[\begin{array}{cc}
T & 0\\
0 & 1\end{array}\right]\ssm{A}{M}{C}{0}\left[\begin{array}{cc}
T^{-1} & 0\\
0 & 1\end{array}\right]=
\sssm{\bar{A}_{11}}{\bar{A}_{12}}{\bar{M}_{1}}%
{\bar{A}_{21}}{\bar{A}_{22}}{\bar{M}_{2}}%
{\bar{C}_{1}}{\bar{C}_{2}}{0}
\]

The \emph{(positive real) balanced truncated model} $\hat{Z}$
is then defined as $\hat{Z}(s)=\bar{C}_1(sI_{m}-\bar{A}_{11})^{-1}\bar{M}_{1}$,
where $m$ is the order of the reduced order system
$\hat{Z}$, i.e. $\bar{A}_{11}\in\Cset^{m\times m}$,
$\bar{M}_{1}\in\Cset^{m\times1}$ and $\bar{C}_{1}\in\Cset^{1\times m}$.

It is important to note that
\[
\left[\begin{array}{cc}
-\bar{A}\Sigma-\Sigma\bar{A}^{*} & \bar{M}-\Sigma\bar{C}^{*}\\
\bar{M}^{*}-\bar{C}\Sigma & 0\end{array}\right]\geq0
\]
and the diagonal structure of $\Sigma$ implies that
\[
\left[\begin{array}{cc}
-\bar{A}_{11}\Sigma_{11}-\Sigma_{11}\bar{A}_{11}^{*} & \bar{M}_{1}-\Sigma_{11}\bar{C}_{1}^{*}\\
\bar{M}_{1}^{*}-\bar{C}_{1}\Sigma_{11} & 0\end{array}\right]\geq0
\]
This ensures that the reduced order model $\hat{Z}$ is
\emph{positive real}, see Lemma~\ref{lem:posreal}!

The order $m$ of the reduced order model may be chosen such that
the approximation error does not exceed an a priori given bound.
In \cite{Green88AC}, equation~(4.30), the following relative error
bound for the spectral densities is given:
\begin{equation}\label{eq:greenbound}
|\Phi(ix)-\hat{\Phi}(ix)|/\Phi(ix) \leq \left(\prod_{k=m+1}^n
\frac{(1+\sigma_k)^2}{(1-\sigma_k)^2}\right) - 1 \mbox{ for all }
x\in \Rset
\end{equation}
where $\Phi=Z+Z^*$ and $\hat\Phi=\hat{Z} + \hat{Z}^*$. Let $k=2c$
be the co-degree of $\Phi(s)$. By the discussion above it follows
that this bound is finite if and only if $m\geq c$ holds.
%Thus one may expect a reasonable approximation only for $m\geq c$.
Furthermore note that for $m \geq c$ the reduced order spectrum
$\hat\Phi$ also has co-degree $2c$, see~\cite{Green88LAA},
Theorem~6.1.

From~(\ref{eq:greenbound}) it is easy to derive an error bound for
the corresponding probability density functions. For simplicity
assume that $\Phi$ is normalized, i.e. $\int \Phi(ix)dx=1$ and
thus $p(x)=\Phi(ix)$ is a pdf. Let
$\hat{p}(x)=\hat{\Phi}(ix)/(\int \hat{\Phi}(ix) dx)$ denote the
approximation of $p(x)$ and let $0<\tau<1$ denote the error bound
on the right hand side of~(\ref{eq:greenbound}).
From~(\ref{eq:greenbound}) we obtain $\Phi(ix)(1-\tau)) \leq
\hat{\Phi}(ix)  \leq  \Phi(ix)(1+\tau))$ and thus
\begin{equation}\label{eq:greenbound2}
(1-\tau) \leq \int_{-\infty}^{\infty} \hat{\Phi}(ix)dx \leq (1+\tau).
\end{equation}
Therefore it follows that
\begin{equation}\label{eq:greenbound3}
|p(x)-\hat{p}(x)|/p(x) \leq \frac{1+\tau}{1-\tau}-1
    = \frac{2\tau}{1-\tau} \mbox{ for all }
x\in \Rset.
\end{equation}
It should be noted that in our experiments we observe
that~(\ref{eq:greenbound2}) is only a rough upper bound
for the ``integrated'' approximation error and thus
(\ref{eq:greenbound3}) is a conservative upper error bound.

In our implementation of the filter a model reduction step is
included after each prediction step. This means after we have
computed a realization of the spectral summand of
$p_{X_{t+1}|Y_1^t}$, we apply the above described scheme to get a
realization of an approximant. This will be used instead of
$p_{X_{t+1}|Y_1^t}$. The order $m$ of the reduced order model is
chosen such that the above error bound~(\ref{eq:greenbound3}) does
not exceed a given threshold $1>\tau >0$. Note that the co-degree
of $p_{X_{t+1}|Y_1^t}$ is bounded by the co-degree of $p_{W_t}$
and that the reduction step does not alter the co-degree!

%%%%%%%%%%%%%%%%%%%%%%%%%%%%%%%%%%%%%%%%%%%%%%%%%%%%%%%%%%%%%%%%%%%%%%%%%%%%
\section{Autocovariance function and estimation\label{sec:ACF}}
%%%%%%%%%%%%%%%%%%%%%%%%%%%%%%%%%%%%%%%%%%%%%%%%%%%%%%%%%%%%%%%%%%%%%%%%%%%%

In this section we analyze the properties of the processes $(Y_t)$
and $(|Y_t|)$. In particular it will be shown,
given some suitable assumptions, that $(Y_t)$ is a white noise
process and that $(|Y_t|)$ is an ARMA process. The mean and the
auto covariance function of $(|Y_t|)$ may be easily computed from
the model parameters in particular from the coefficients of the
polynomial $V(x)$ and from the moments of the noise processes
$(W_t)$ and $(U_t)$. This enables the use a simple method of
moments to estimate the model parameters.

The standing assumptions in this section are as follows:
\begin{enumerate}
\item $V(x)=v_0+v_1 x + \cdots v_d x^d$ is a non negative polynomial
    ($V(x)\geq 0$ for all $x \in \Rset$) and it has order
    $d$.
\item The processes $(W_t)$ and $(U_t)$ are two i.i.d processes, which
    are independent from each other. The moments $M_W(k):=\Ex W_t^k$ exist
    for all $0\leq k \leq m_W$ and $m_W\geq 2d$ holds. The moments
    $M_U(k):=\Ex U_t^k$ exist for all $0 \leq k \leq m_U$ and
    $m_U\geq 2$ holds.
\item The parameter $a$ is bounded by $|a|<1$.
\end{enumerate}
Note that within this section it is not needed that $W_t$ and $U_t$ have
rational probability density functions.

The main result is given in the following Proposition:

\begin{proposition}\label{prop:ACF} Under the assumptions (i),
(ii) and (iii) there  exists a strictly
stationary solution \( (X_{t},Y_t) \) of the model
(\ref{eq:model}).

The moments \( M_X(k):=\Ex X^{j}_{t} \) exist up to order \( m_{X}=m_{W} \)
and may be computed recursively from the relations (starting with
\( M_X(0)=1 \))
\begin{equation}\label{eq:momX}
M_X(k)=\frac{1}{1-a^{k}}\sum ^{k-1}_{l=0}\binom{k}{l}a^{l}M_W(k-l)M_X(l)\; \: ;\; \: 1\leq k\leq m_{W}
\end{equation}

The process \( (Y_{t}) \) is a white noise process.

The process \( Z_{t}=|Y_{t}| \) is an ARMA process of order less
than or equal to \( d+1 \).
\end{proposition}

\proof Let \( M_{t,k}=\sum ^{k}_{j=1}|a|^{j-1}|W_{t-j}| \). Since
\( M_{t,k} \) is monotonically increasing with \( k \), and since
\( \Ex M_{t,k}\leq \Ex |W_{t}|/(1-|a|) \) is bounded, we conclude
that \( \lim _{k\rightarrow \infty }M_{t,k} \) and \( X_{t}:=\lim
_{k\rightarrow \infty }\sum ^{k}_{j=1}a^{j-1}W_{t-k} \) exist a.s.
Furthermore \( \Ex X_{t}=\Ex W_{t}/(1-a) \).

Now suppose that \( \Ex X^{l}_{t} \) exists for \( 1\leq l<k\leq
m_{W} \).  From (\ref{eq:model}) it follows that
\[
X^{k}_{t+1}-a^{k}X^{k}_{t}=
\sum ^{k-1}_{l=0}\binom{k}{l}a^{l}X^{l}_{t}W^{k-l}_{t}.
\]
Since \( \sum ^{k-1}_{l=0}\binom{k}{l}a^{l}X^{l}_{t}W^{k-l}_{t} \)
has a finite mean it follows analogously that \( \Ex X^{k}_{t} \)
exists. By taking expectations on both sides of the above equation
and by using the independence of \( X_{t} \) and \( W_{t} \) one
obtains (\ref{eq:momX}).

Furthermore for \( k\geq 0 \), \( X^{i}_{t+k+1}X^{j}_{t}=(aX_{t+k}+W_{t+k})^{i}X^{j}_{t}=
\sum ^{i}_{l=0}\binom{i}{l}a^{l}W^{i-l}_{t+k}X^{l}_{t+k}X^{j}_{t} \)
and thus \( \Ex (X^{i}_{t+k+1}X^{j}_{t})=
\sum ^{i}_{l=0}\binom{i}{l}a^{l}M_W(i-l)\Ex (X^{l}_{t+k}X^{j}_{t}) \).

Define \( \vec{X}_{t}=(1,X_{t},\ldots \, X^{d}_{t})' \),
\( \vec{M}_X = \Ex \vec{X}_{t}=(1,M_X(1),\ldots,M_X(d))' \),
\( \vec{V}=(v_{0},\ldots \, v_{d})' \) and
\[
F=\left( \begin{array}{ccccc}
\binom{0}{0}a^{0}M_W(0) & 0 & \cdots  & \cdots  & 0\\
\binom{1}{0}a^{0}M_W(1) & \binom{1}{1}aM_W(0) & \ddots  &  & \vdots \\
\binom{2}{0}a^{0}M_W(2) & \binom{2}{1}aM_W(1) & \binom{2}{2}a^{2}M_W(0) & \ddots  & \vdots \\
\vdots  &  &  & \ddots  & 0\\
\binom{d}{0}a^{0}M_W(d) & \cdots  & \cdots  & \cdots  & \binom{d}{d}a^{d}M_W(0)
\end{array}\right)
\]

Using these notations the above relations may be written as:
\( \vec{M}_X=F\vec{M}_X \) and
\( \Ex \vec{X}_{t+k+1}\vec{X}_{t}'=F\Ex \vec{X}_{t+k}\vec{X}_{t}' \).

First consider the process \( (V(X_{t})) \). It is immediate to
see that \( \Ex V(X_{t})=\vec{V}'\vec{M}_X \) and \( \Ex
V(X_{t+k})V(X_{t})=\vec{V}'F^{k}(\Ex
\vec{X}_{t}\vec{X}_{t}')\vec{V} \), for \( k\geq 0 \). Note that
\( F \) has eigenvalues \( 1,a,\ldots,a^{d} \) and that \(
e=(1,0,\ldots ,0) \) and \( \vec{M}_X \) are the left and the
right eigenvectors corresponding to the eigenvalue \( 1 \).
Furthermore \( e\Ex \vec{X}_{t}\vec{X}_{t}'=\vec{M}_X' \). This
implies that the auto-covariance function of \( V(X_{t}) \) is
given by
\[
\mathrm{Cov}(V(X_{t+k}),V(X_{t}))=\left\{
\begin{array}{ll}
\vec{V}'(\Ex \vec{X}_{t}\vec{X}_{t}'-\vec{M}_X\vec{M}_X')\vec{V} &
\textrm{for}\, k=0\\
\vec{V}'(F-\vec{M}_X e)^{k-1}
    \left((F-\vec{M}_X e)\Ex \vec{X}_{t}\vec{X}_{t}'\vec{V}\right) &
        \textrm{for}\, k>0
\end{array}
\right.
\]
From the above representation it follows that \( (V(X_t)) \) is an
ARMA process of order less than or equal to \( d+1 \). Note that
\( (F-\vec{M}_X e) \) has eigenvalues \( 0,a,\ldots,a^{d} \).

Next consider the process \( (Y_{t}). \) We have
\( \Ex Y_{t}=\Ex V(X_{t})\Ex U_{t}=0 \),
by the independence of \( X_{t} \) and of \( U_{t} \). The auto
covariance function of \( (Y_{t}) \) is given by
\[
\Ex Y_{t+k}Y_{t}=\Ex (V(X_{t+k})V(X_{t}))\Ex (U_{t+k}U_{t})=
\left\{ \begin{array}{ll}
\Ex V(X_{t})^{2}\Ex U^{2}_{t} & \textrm{for}\, k=0\\
0 & \textrm{for}\, k>0
\end{array}\right.
\]

Finally let us consider \( |Y_{t}| \). The mean value of \( |Y_t| \)
is \( \Ex |Y_{t}|=\Ex V(X_{t})\Ex |U_{t}| \) and the second moments
are given by \( \Ex |Y_{t}|^{2}=\Ex V(X_{t})^{2}\Ex U^{2}_{t} \)
and \( \Ex |Y_{t+k}||Y_{t}|=\Ex (V(X_{t+k})V(X_{t}))(\Ex |U_{t}|)^{2} \)
for \( k>0 \). This implies
\[
\mathrm{Cov}(|Y_{t+k}|,|Y_{t}|)=\left\{ \begin{array}{ll}
\Ex V(X_{t})^{2}\Ex U^{2}_{t}-(\Ex V(X_{t})\Ex |U_{t}|)^{2} &
    \textrm{for}\, k=0\\
\mathrm{Cov}(V(X_{t+k}),V(X_{t}))(\Ex |U_{t}|)^{2} &
    \textrm{for}\, k>0
\end{array}\right.
\]
\eproof

Of course analogous calculations apply for \( |Y_{t}|^{k} \),
provided, that sufficiently many moments of \( W_{t} \) and of \( U_{t} \)
exist.

%%%%%%%%%%%%%%%%%%%%%%%%%%%%%%%%%%%%%%%%%%%%%%%%%%%%%%%%%%%%%%%%%%%%%%%%%%%%%
\section{Simulation results}
%%%%%%%%%%%%%%%%%%%%%%%%%%%%%%%%%%%%%%%%%%%%%%%%%%%%%%%%%%%%%%%%%%%%%%%%%%%%%

All simulation and estimation results presented here are based on the
following specifications:

\begin{equation}
\label{eq:simmodel}
\begin{array}{rcl}
X_{t+1} & = & aX_{t}+ W_{t}\\
Y_{t} & = & \Psi V(\sigma X_{t})U_{t}
\end{array}
\end{equation}
This is a slight reformulation of the model~(\ref{eq:model}). The idea is
to fix the function $V(x)$ and the distributions of $W_t$ and of $U_t$,
which leaves the three parameters $a,\Psi$ and $\sigma$ for estimation.

The function $V(x)$ is chosen as
\[
V(x)=(1+\frac{x}{2d})^{d}+0.1
\]
which is a rough approximation of $\exp(x/2)$. The additional constant
$0.1$ is added to ensure $V(x)>0$.

The inputs $W_t$ and $U_t$ are assumed to have scaled
t-distributions. This means that $c_W W_t$ has a t-distribution
with $n_W$ degrees of freedom and the scaling constant chosen such
that  $\Ex W_t^2=1$ holds. Analogously $c_U U_t$ has a
t-distribution with $n_U$ degrees of freedom and the scaling is
such that $\Ex U_t^2=1$ holds. Throughout this section the integer
parameters $d, n_W$ and $n_U$ are fixed and given by $d=4$,
$n_W=9$ and $n_U=3$. This implies in particular that the
assumptions of Proposition~\ref{prop:ACF} are fulfilled.

The first part of this section deals with the estimation of the
parameters $(a,\Psi,\sigma).$ Table~\ref{tab:moments} shows the moments of
the process \( (|Y_{t}|) \) for some combinations of the
parameters \( a,\, \Psi=1.0 ,\, \sigma  \).

In a small simulation study we have investigated the performance
of a simple method-of-moments estimation, where \( 10 \) lags of
the auto-covariance function of \( |Y_{t}| \) have been used. To
be more precise let
\[
m(a,\Psi,\sigma): = (\Ex |Y_t|,\mathrm{Var} |Y_t|,
    \mathrm{Cov}(|Y_{t+1}|,|Y_{t}|),\ldots,\mathrm{Cov}(|Y_{t+10}|,|Y_{t}|))
\]
and let $\hat{m}_T$ be the sample estimate of this vector of moments
given a sample of size $T$. Then the estimates $(\hat{a},\hat{\Psi},\hat{\sigma})$
are computed by minimizing
\[
\| m(a,\Psi,\sigma)-\hat{m}_T \|^2
\]
The results for \( 1000 \) simulation runs for
(simulated) data series of length \( T=1000 \) are collected in table~\ref{tab:simresults}.
Both the mean estimation error (mean) and the standard deviation (std) over
these 1000 simulation are shown in dependence of the true parameters.
Note e.g. that the estimate of \( a \) shows a significant bias especially
for small \( a \) and \( \sigma  \). However this is only a first
rough estimation scheme and other enhanced estimates will be investigated
in future.

It has been mentioned in section~\ref{sec:filter} that the filter
is able to compute the likelihood. However the computation of the
filter is presently too time demanding to
implement a maximum likelihood estimation based on the filter.

Next we test the filter on some real world data. In particular we
consider data which also have been analysed
by~\cite{MahieuSchotman98}. The authors consider five exchange
rate data series and study the empirical performance of stochastic
volatility models. Here we only consider the Dollar/Yen exchange
rate data, which consists of $T=1102$  weekly observations from 3
January 1973 until 9 February 1994.

The parameters of the model~(\ref{eq:simmodel}) have been
estimated by the method of moments as described in the previous
section. Here 25 lags of the auto covariance are used and the
obtained estimates are $\hat{a} = 0.957$, $\hat{\sigma} =
0.309$ and $\hat{\Psi} = 0.921$. Figure~\ref{fig:yen_dollar2}
shows the sample ACF and the fitted ACF. Next the filter is run on
this data set to compute a one step ahead prediction of
$|Y_{t+1}|$. The result is shown in figure~\ref{fig:yen_dollar3}.
Note that, since the stationary solution of the state $X_t$ is not
rationally distributed as far as we know, we have simply assumed
that $X_1$ has a scaled t-distribution with $n_X=9$ degrees of
freedom and the scaling was chosen such that the variance of $X_1$
is equal to $1/(1-a^2)$ i.e. equal to variance of the stationary
solution.

Finally we consider some simulated data. The parameters
were chosen as $a=0.9$, $\Psi=2$ and $\sigma=1.5$. The simulation and
the filter were initialized with a scaled t-distributed
random variable $X_1$, where the degrees of freedom
is $n_X=9$ and the scaling parameter is chosen such that
the variance of $X_1$ is equal to $1/(1-a^2)$. The length of the
simulated series is $T=100$.

Figure~\ref{fig:simulatedData_fig3} shows the simulated trajectory of $Y_t$ and the one step ahead prediction of $|Y_t|$, i.e.
$\Ex (|Y_t|\,|\,Y_1^{t-1})=
    \Psi (\Ex |U_t| ) \Ex( V(\sigma X_t) \,|\,Y_1^{t-1} )$.
The conditional expectation $\Ex (V(\sigma X_t) \,|\,Y_1^{t-1})$
is computed from the conditional probability density
function $p_{X_t|Y_1^{t-1}}$, which is computed by the
filter. See also Propositions~\ref{prop:Ex} and~\ref{prop:ACF}.

Figure~\ref{fig:simulatedData_fig6} shows the conditional probability density function $p_{X_{t+1}|Y_1^{t}}$, for $t=100$.
In each time step of the filter balanced model reduction
is used as described above. Let $\hat{p}_{X_{t+1}|Y_1^{t}}$ denote the
approximation of the conditional pdf $p_{X_{t+1}|Y_1^{t}}$.
The order $m$ of the the reduced order system is chosen such that the
relative error $|p_{X_{t+1}|Y_1^{t}}(x) - \hat{p}_{X_{t+1}|Y_1^{t}} | / p_{X_{t+1}|Y_1^{t}}(x)$
is at most $0.02$, i.e. we allow at most an error of 2 percent.
See equation~(\ref{eq:greenbound3}).
For this specific example typical model orders are $n=85$ and $m=9$,
which means that the state dimension is almost reduced by a factor $10$.
If one compares the conditional expectation of $X_t$ and of $V(X_t)$ given the
observations $Y_1^t$ computed from the full order pdf $p_{X_{t+1}|Y_1^{t}}$ and
from the approximant $\hat{p}_{X_{t+1}|Y_1^{t}}$ then in this example the relative error is of the order $10^{-14}$. These numbers
indicate the excellent quality of the used approximation scheme.

Finally figure~\ref{fig:simulatedData_fig7} shows the
evolution of the conditional densities
$p_{X_{t+1}|Y_1^{t}}$ over time.

%%%%%%%%%%%%%%%%%%%%%%%%%%%%%%%%%%%%%%%%%%%%%%%%%%%%%%%%%%%%%%%%%%%%%%%%%%%%%%
\section{Conclusion}
%%%%%%%%%%%%%%%%%%%%%%%%%%%%%%%%%%%%%%%%%%%%%%%%%%%%%%%%%%%%%%%%%%%%%%%%%%%%%%

The exact filter for a class of stochastic volatility models is
derived. A standard stochastic volatility model in which the
disturbances are Gaussian and the volatility function involved is
exponential can be viewed as a limiting case.  The complexity of
the exact filter increases in the sense that the matrices that are
used to represent the rational probability density functions  tend
to grow quickly. An approximate filter is presented in which at
each time step the conditional probability density function of the
state, which is rational, is replaced by an approximating rational
probability density function, using the SBT method (stochastically
balanced truncation). Using a well-known error bound the
approximating rational probability density function can be chosen
such that on each point of the real line the relative error is
less than a given percentage (the tolerance level involved can be
chosen by the user). In some simulated and empirical applications
we find that using a tolerance level of as low as 2 percent still
leads to an enormous reduction in complexity, keeping the order of
the rational functions well within bounds that are considered
tractable with modern computers. Lower tolerance levels could also
be achieved if desired, but then larger matrices will have to be
handled. The model presented is very flexible, especially with
respect to the specification of the probability density functions
for the disturbances. Here one can vary between very heavy-tailed
disturbances (with Cauchy density for instance) and less
heavy-tailed disturbances (with Student-t densities that are
approximating Gaussian densities for example). In the applications
in this paper we have stayed as close as possible to the
traditional Gaussian model. However the possibility of specifying
more heavy-tailed densities seems one of the most interesting
features of this class of models. Exploring those possibilities is
an interesting topic for future research. Also valuation of
financial derivatives in a market in which the asset price
movements can be described by a stochastic volatility model of the
type investigated here, is an interesting topic for future
research. More generally the methodology of working with rational
density functions in filtering problems in the way presented here
could have a much wider range of applications, as the methodology
is really general and flexible and numerically stable methods for
various operations involved are now provided.  Preliminary
experience with the methodology shows especially striking results
deriving from the application of the SBT approximation method. It
is to be expected that this can also be successfully applied to
the linear filtering problems with rationally distributed
disturbances considered in \cite{HanzonOber01}.

%%%%%%%%%%%%%%%%%%%%%%%%%%%%%%%%%%%%%%%%%%%%%%%%%%%%%%%%%%%%%%%%%%%%%%%%%%%%
\appendix
\section{Results from system theory}
%%%%%%%%%%%%%%%%%%%%%%%%%%%%%%%%%%%%%%%%%%%%%%%%%%%%%%%%%%%%%%%%%%%%%%%%%%%%

%%%%%%%%%%%%%%%%%%%%%%%%%%%%%%%%%%%%%%%%%%%%%%%%%%%%%%%%%%%%%%%%%%%%%%%%%%%%
\subsection{Numerical calculation of the co-degree and of the zeros
of a strictly proper rational function\label{app:codegree}}
%%%%%%%%%%%%%%%%%%%%%%%%%%%%%%%%%%%%%%%%%%%%%%%%%%%%%%%%%%%%%%%%%%%%%%%%%%%%

Consider a strictly proper scalar rational function%
\footnote{In this section $G$ is an arbitrary, not necessarily stable,
transfer function. We will use results obtained in this section e.g.
for a spectrum $\Phi=KK^{*}$ and for its factor $K^{*}$.}
\[
G(s)=C(sI_{n}-A)^{-1}B=\frac{a_{0}+a_{1}s+\cdots + a_{q}s^{q}}%
{b_{0}+b_{1}s+\cdots+ b_{n}s^{n}}
\]
where
\[
q<n\mbox{ and }
a(s)=a_{0}+a_{1}s+\cdots + a_{q}s^{q},\;
b(s)=b_{0}+b_{1}s+\cdots + b_{n}s^{n}\;
\mbox{ are coprime}.
\]

The \emph{co-degree} of $G$ is defined as $(n-q)$, i.e. as the
multiplicity of the infinite zero of $G(s)$. Since $G$ is strictly
proper the co-degree is positive. The Taylor series expansion of
$G(s)$ at infinity is given by
\begin{eqnarray*}
G(s) & = & G_{0}+G_{1}s^{-1}+G_{2}s^{-2}+\cdots\\
 & = & D+CBs^{-1}+CABs^{-2}+\cdots
 \end{eqnarray*}
Therefore the co-degree of $G$ is related to the Markov parameters
of $G$ as follows.

\begin{lemma}\label{lem:codegG}
The co-degree of $G(s)=C(sI_{n}-A)^{-1}B$ is equal to $c$ iff
\begin{itemize}
\item $CA^{c-1}B\neq0$ and $CA^{i-1}B=0$ for all $1\leq i<c$.
\end{itemize}
\end{lemma}

Note that a naive check on $CA^{i-1}B=0$ in order to compute the
co-degree is numerically unstable since $A$ might have eigenvalues
of modulus larger than one and thus round off errors would
``explode''.

The (finite) zeros of the transfer function $G(s)$ are the (finite)
eigenvalues of the pencil:
\begin{equation}
\lambda E-N :=
\left[\begin{array}{cc}
\lambda I_{n}-A & B\\
-C & D\end{array}\right]\label{eq:pencilG}
\end{equation}

Therefore the co-degree and the finite zeros of $G$ may be
computed from the eigenstructure of the above pencil. We will make
use of the following concepts, see e.g. \cite{VanDooren81}. A
pencil $(\lambda E-N)$ is called regular if it is square and if
$\det(\lambda E-N)$ is not constant. The zeros of $\det(\lambda
E-N)$ are the eigenvalues of the pencil. Suppose there exist full
column rank matrices $X, Y \in \Cset^{n\times k}$, $k\leq n$ and
matrices $\bar{E},\bar{N}\in \Cset ^{k\times k}$ such that
\begin{equation}
(\lambda E-N)X=Y(\lambda \bar{E}-\bar{N})
\end{equation}
holds. The space spanned by the columns of $X$ is called a deflating subspace of
the pencil $(\lambda E-N)$. This is a generalization of the concept
of invariant subspaces to arbitrary pencils. The $k$-dimensional
pencil $(\lambda \bar{E}-\bar{N})$ is called a divisor of $(\lambda E-N)$.
If $\bar{E}$ is non singular, then $(\lambda \bar{E}-\bar{N})$ is called
a finite divisor of $(\lambda E-N)$. In this case $(\lambda \bar{E}-\bar{N})$
has $k$ finite eigenvalues which are of course also eigenvalues of $(\lambda E-N)$.
If there exist two  non singular matrices $S,T \in \Cset^{k\times k}$ such that
\[
S(\lambda \bar{E}-\bar{N})T=
\left(\begin{array}{ccccc}
    -1 & \lambda &       0 & \cdots &      0 \\
    0  &      -1 & \lambda & \ddots & \vdots  \\
\vdots &  \ddots &  \ddots & \ddots & \vdots \\
    0  &  \cdots &       0 &     -1 & \lambda \\
    0  &  \cdots &  \cdots &      0 &   -1 \\
\end{array}\right)
\]
then $(\lambda \bar{E}-\bar{N})$ is called an elementary infinite divisor.

An alternative characterisation of the co-degree now is as
follows:

\begin{lemma}
The co-degree of $G(s)$ is positive and it is equal to $c$ iff the
pencil~(\ref{eq:pencilG}) has an elementary infinite divisor of
dimension $(c+1)$ and a finite divisor of dimension $(n-c)$.
\end{lemma}

\proof To prove this lemma the pencil is transformed to a socalled
\emph{staircase} form as defined in \cite{VanDooren81}. This will
also give a numerically robust way to analyze the co-degree and
the eigenstructure of the above pencil.

Let $U_{1}\in\Cset^{n\times n}$ be a row compression of $(-B)$,
i.e. $U_{1}$ is a unitary matrix such that $U_{1}^{*}(-B)=[\bar{b},0,\ldots,0]^{*}$
and $\bar{b}>0$. (Note that $B\neq0$.) Apply this state space transformation
and define \[
\left[\begin{array}{cc}
U_{1}^{*} & 0\\
0 & 1\end{array}\right]\ssm{A}{-B}{C}{-0}\left[\begin{array}{cc}
U_{1} & 0\\
0 & 1\end{array}\right]=:\ssm{A_{1}}{-B_{1}}{C_{1}}{0}\]

Note that $CB=C_{1}B_{1}$ and thus the first element of $C_{1}$
is zero iff $c>1$. In the next step let \[
U_{2}=\left[\begin{array}{cc}
1 & 0\\
0 & \bar{U}_{2}\end{array}\right]\]
where $\bar{U}_{2}\in\Cset^{n-1\times n-1}$ is a row compression
of the last $n-1$ entries of the first column of $A_{1}$. Apply
this state space transformation to get \[
\left[\begin{array}{cc}
U_{2}^{*} & 0\\
0 &
1\end{array}\right]\ssm{A_{1}}{-B_{1}}{C_{1}}{0}\left[\begin{array}{cc}
U_{2} & 0\\
0 & 1\end{array}\right]=:\ssm{A_{2}}{-B_{2}}{C_{2}}{0}\]

By construction the $(1,1)$ element of $A_{1}$ and the first elements
of $C_{1}$ and of $B_{1}$ are not affected by this transformation.
Furthermore note that the last $n-1$ elements of $B_{2}$ and the
last $n-2$ elements of the first row of $A_{2}$ are zero. In addition
we have $CAB=C_{2}A_{2}B_{2}=0$ iff $c>2$.  Thus $c>2$ holds iff the second
element of $C_2$ is zero.

Now this procedure is repeated until a nonzero element pops up in
the $k$-th position of $C_k$. This is a possible way to estimate
the co-degree of $G$.

After $c+1$ steps of this kind we end up with a matrix of the form:

\begin{equation}
\left[\begin{array}{cccccc|cccccc|c}
* & * & \cdots & \cdots & \cdots & * & * & \cdots &  &  & \cdots & * & \oplus\\
\oplus & * & \cdots & \cdots & \cdots & * & * & \cdots &  &  & \cdots & * & 0\\
0 & \ddots & \ddots &  &  & \vdots & \vdots &  &  &  &  & \vdots & \vdots\\
\vdots & \ddots & \ddots & \ddots & & \vdots & \vdots & & & & & \vdots &  \vdots \\
\vdots &  & \ddots & \ddots & \ddots & \vdots & \vdots &  &  &  &  & \vdots & \vdots\\
0 & \cdots & \cdots & 0 & \oplus & * & * & \cdots &  &  & \cdots & * & 0\\
\hline
0 & \cdots & \cdots & \cdots & 0 & \beta & * & \cdots &  &  & \cdots & * & 0\\
0 & \cdots & \cdots & \cdots & 0 & 0     & * & \cdots &  &  & \cdots & * & \vdots\\
\vdots &  &  &  & \vdots & \vdots & \vdots &  &  &  &  & \vdots & \vdots\\
0 & \cdots & \cdots & \cdots & 0 & 0 & * & \cdots &  &  & \cdots & * & 0\\
\hline
0 & \cdots & \cdots & \cdots & 0 & \alpha & * & \cdots &  &  & \cdots & * & 0\end{array}\right]
\label{eq:staircase}\end{equation}

The horizontal and vertical lines partition the above matrix into
blocks of size $c$, $n-c$ and $1$ respectively. Two particular
elements of the above matrix, namely the $(c+1,1)$ and the $(n+1,c)$
element, are denoted with $\beta$ and $\alpha$ respectively. Note
that $\beta>0$ and $\alpha\neq0$ holds.

Note that (for $j<c$) the columns $[1,\ldots,j+1]$ of the matrix
$U=U_{1}U_{2}\cdots \cdots U_{c+1}$ form an orthogonal basis of the
column space of $[B,AB,\ldots,A^{j}B]$.

By a permutation of rows and columns we bring the last column to the
first position and the last row to the $(c+1)$-th position. Finally
apply the Givens rotation \[
\bar{Q}=\left[\begin{array}{cc}
\alpha^{*} & \beta^{*}\\
-\beta & \alpha\end{array}\right]\frac{1}{\sqrt{\alpha^{*}\alpha+\beta^{*}\beta}}=:\left[\begin{array}{cc}
q_{11} & q_{12}\\
q_{21} & q_{22}\end{array}\right]\]
to the rows $c+1$ and $c+2$. If $Q$ and $Z$ denote the concatenation
of all these unitary row and column operations, then we have

\begin{equation}\label{eq:staircase2}
Q\ssm{A}{-B}{C}{0}Z= \left[\begin{array}{ccccc|ccc}
\oplus & *       & \cdots & \cdots & *      & *      & \cdots &  *\\
0      & \ddots  & \ddots &        & \vdots & \vdots &        & \vdots\\
\vdots & \ddots  & \ddots & \ddots & \vdots & \vdots &        & \vdots\\
\vdots &         & \ddots & \ddots & *      & \vdots &        & \vdots\\
0      & \cdots  & \cdots & 0      & \oplus & *      & \cdots &  *\\
\hline
0      & \cdots  & \cdots & \cdots & 0      & *      & \cdots & *\\
\vdots &         &        &        & \vdots & \vdots &        & \vdots\\
0      & \cdots  & \cdots & \cdots & 0      & *      & \cdots  & *
\end{array}\right]=:
\left[\begin{array}{cc}
\bar{N}_{11} & \bar{N}_{12}\\
           0 & \bar{N}_{22}
\end{array}\right]
\end{equation}

and

\begin{equation}\label{eq:staircase3}
Q\ssm{I_{n}}{0}{0}{0}Z=
\left[\begin{array}{ccccc|ccccc}
0     &     1&0     &\cdots&0     &0     &0     &\cdots&\cdots&0\\
0     &\ddots&\ddots&\ddots&\vdots&\vdots&\vdots&      &      &\vdots \\
\vdots&\ddots&\ddots&\ddots&0     &\vdots&\vdots&      &      &\vdots \\
\vdots&      &\ddots&\ddots&1     &0     &0     &\cdots&\cdots&\vdots \\
0     &\cdots&\cdots&0     &0     &q_{12}&0     &\cdots&\cdots&0\\
\hline
0     &\cdots&\cdots&\cdots&0     &q_{22}&0     &\cdots&\cdots&0\\
\vdots&      &      &      &\vdots&0     &1     &\ddots&      &0\\
\vdots&      &      &      &\vdots&\vdots&\ddots&\ddots&\ddots&\vdots\\
\vdots&      &      &      &\vdots&\vdots&      &\ddots&\ddots&0\\
0     &\cdots&\cdots&\cdots&0     &0     &\cdots&\cdots&0     &1
\end{array}\right]=:\left[\begin{array}{cc}
\bar{E}_{11} & \bar{E}_{12}\\
           0 & \bar{E}_{22}
\end{array}\right]
\end{equation}

Now this block upper triangular form displays the eigenstructure
of the pencil $(\lambda E-N)$. Since
$\bar{N}_{11}\in\Cset^{c+1\times c+1}$ is an upper-triangular
non-singular matrix it follows that
$(\lambda\bar{E}_{11}-\bar{N}_{11})$ is an $(c+1)$ dimensional
elementary infinite divisor of the pencil. Furthermore note that
$Z$ may be partitioned as
\[
Z = \left[\begin{array}{cc} 0 & U \\ 1 & 0 \end{array}\right]
\]
and that the first $c+1$ columns of $Z$ form a basis for the
deflating subspace corresponding to this infinite divisor. The
same holds true if we only take the first $j+1$ columns, for
$0\leq j\leq c$. To be more precise consider the $(j+1\times j+1)$
dimensional left upper sub-block of
$(\lambda\bar{E}_{11}-\bar{N}_{11})$. By the triangular structure
of the matrices $\bar{E}_{11}$ and $\bar{N}_{11}$ it follows that
this sub-block defines an infinite elementary divisor and that the
first $j+1$ columns of $Z$ span the corresponding deflating
subspace.

Since $\bar{E}_{22}\in\Cset^{(n-c)\times (n-c)}$ is non singular
it follows that $(\lambda\bar{E}_{22}-\bar{N}_{22})$ is an $(n-c)$
dimensional finite divisor of the pencil. \eproof

%%%%%%%%%%%%%%%%%%%%%%%%%%%%%%%%%%%%%%%%%%%%%%%%%%%%%%%%%%%%%%%%%%%%%%%%%%%
\subsection{Elementary operations on rational functions\label{app:add+mult}}
%%%%%%%%%%%%%%%%%%%%%%%%%%%%%%%%%%%%%%%%%%%%%%%%%%%%%%%%%%%%%%%%%%%%%%%%%%%

Let two strictly proper rational function
$G_i=C_i(sI-A_i)^{-1}B_i$ with state space realizations
$(A_i,B_i,C_i)$, $i=1,2$ be given.

A state space realization for $G_1^*(s)=B_1^*(-sI-A_1^*)^{-1}C_1^*$ is
given by
\[
G_1^*=\pi \ssm{-A_1^*}{C_1^*}{-B_1^*}{0}
\]
The sum $G_1+G_2$ has a state space realization:
\[
G_1+G_2=\pi \sssm{A_1}{0}{B_1}{0}{A_2}{B_2}{C_1}{C_2}{0}
\]
The product $G_1G_2$ has a state space realization:
\[
G_1G_2=\pi \sssm{A_1}{B_1B_2}{0}{0}{A_2}{B_2}{C_1}{0}{0}
\]
If $C_1B_1=0$ then $G_1(ys)s$ is strictly proper and a state space realization
is given by
\[
G_1(ys)s=\pi \ssm{A_1y^{-1}}{A_1B_1y^{-1}}{B_1y^{-1}}{0}
\]

%%%%%%%%%%%%%%%%%%%%%%%%%%%%%%%%%%%%%%%%%%%%%%%%%%%%%%%%%%%%%%%%%%%%%%%%%%%%
\subsection{Computation of a spectral summand from a spectral factor
\label{app:rfac2rsum}}
%%%%%%%%%%%%%%%%%%%%%%%%%%%%%%%%%%%%%%%%%%%%%%%%%%%%%%%%%%%%%%%%%%%%%%%%%%%%

Suppose we have given  a (stable) spectral factor $K(s)=C(sI_{n}-A)^{-1}B$
and that we want to compute a spectral summand of $\Phi(s)=K(s)K^{*}(s)$:

Let $P$ be the solution of the Lyapunov equation
\[
AP+PA^{*}+BB^{*}=0
\]
and define $M=PC^{*}$. The state space transformation $T$
\begin{equation}
T=\left[\begin{array}{cc}
I_{n} & P\\
0 & I_{n}\end{array}\right]\label{eq:T}
\end{equation}
then gives
\[
\left[\begin{array}{cc}
I_{n} & P\\
0 & I_{n}
\end{array}\right]\left[\begin{array}{cc}
A & -BB^{*}\\
0 & -A^{*}
\end{array}\right]\left[\begin{array}{cc}
I_{n} & -P\\
0 & I_{n}
\end{array}\right]=\left[\begin{array}{cc}
A & -AP-PA^{*}-BB^{*}\\
0 & -A^{*}
\end{array}\right]=\left[\begin{array}{cc}
A & 0\\
0 & -A^{*}
\end{array}\right],\]

\[
\left[\begin{array}{cc}
I_{n} & P\\
0 & I_{n}
\end{array}\right]\left[\begin{array}{c}
0\\
C^{*}
\end{array}\right]=\left[\begin{array}{c}
PC^{*}\\
C^{*}
\end{array}\right]=\left[\begin{array}{c}
M\\
C^{*}\end{array}\right], \]
\[
\left[\begin{array}{cc}
C & 0
\end{array}\right]\left[\begin{array}{cc}
I_{n} & -P\\
0 & I_{n}
\end{array}\right]=\left[\begin{array}{cc}
C & -CP
\end{array}\right]=\left[\begin{array}{cc}
C & -M^{*}\end{array}\right], \] and thus $Z(s)=C(sI_{n}-A)^{-1}M$
is a (stable) spectral summand of $\Phi(s)$.

%%%%%%%%%%%%%%%%%%%%%%%%%%%%%%%%%%%%%%%%%%%%%%%%%%%%%%%%%%%%%%%%%%%%%%%%
\subsection{Computation of a spectral summand from an
integrable spectral density\label{app:rpdf2rsum}}
%%%%%%%%%%%%%%%%%%%%%%%%%%%%%%%%%%%%%%%%%%%%%%%%%%%%%%%%%%%%%%%%%%%%%%%%%%%

Let $\Phi=H(sI_{2n}-F)^{-1}G$ be given. First compute a Schur decomposition
of $F$ such that the stable eigenvalues of $F$ appear on the first $n$
positions, i.e.
\[
\bar{F}=V^* F V=
\left[\begin{array}{cc}
\bar{F}_{11} & \bar{F}_{12} \\
0 & \bar{F}_{22}
\end{array}\right]
\]
where $V$ is a unitary matrix, $\bar{F}$ is an upper triangular matrix
and $\bar{F}_{11}\in\Cset^{n\times n}$ is asymptotically stable.

Solve the Lyapunov equation
\[
-\bar{F}_{11}P+P\bar{F}_{22}+\bar{F}_{12}=0
\]
and set
\[
\begin{array}{rcl}
A &=& \bar{F}_{11}\\
M &=& [I,P]V^* G \\
C &=& HV (I_n,0)^*
\end{array}
\]
to get a stable spectral summand $Z=C(sI-A)^{-1}M$.

%%%%%%%%%%%%%%%%%%%%%%%%%%%%%%%%%%%%%%%%%%%%%%%%%%%%%%%%%%%%%%%%%%%%%%%%%%%%%
\subsection{Operations on rational densities\label{app:ratdensop}}
%%%%%%%%%%%%%%%%%%%%%%%%%%%%%%%%%%%%%%%%%%%%%%%%%%%%%%%%%%%%%%%%%%%%%%%%%%%%%

In \cite{HanzonOber01} it was show that the operations of
translation, scaling, multiplication and convolution of rational
densities can be translated into linear algebra operations on
corresponding state-space realizations of spectral summands. For
ease of reference here we give some of these results which are
needed for the implementation of the filter. Note that
multiplication of two rational functions could be implemented via
their summands. However using spectral factors seems to be
numerically more reliable. Thus in our implementation of the
filter we have chosen this approach.

\begin{proposition}\label{prop:scaling+convolution}
Let $X_1$ and $X_2$ denote stochastically independent random
variables with rational density functions
 $p_1,~p_2.$
 For $j=1,2,$ let $Z_j(s)$ denote the corresponding stable spectral summand,
with a state-space realization $[A_j,M_j,C_j]$ with state-space
dimension $n_j.$
\begin{itemize}
\item[(i)] For $a\neq 0$ the random variable  $X=aX_1$ has a
rational density whose spectral summand has a state space
realization given by $[A,M,C]=[aA_1,M_1,C_1]$ if $a>0$ and
$[A,M,C]=[-aA_1^*,C_1^*,M_1^*]$ if $a<0.$ \item[(ii)] The sum
$X=X_1+X_2$ has a rational density function $p=p_1 \star p_2$,
i.e. the convolution of $p_1$ and $p_2$, and the spectral summand
of $p$ has a state-space realization given by $[A,M,C]$ where
 $A=A_1 \otimes I_{n_2}+I_{n_1} \otimes A_2,$
 $M=M_1 \otimes M_2$ and $C=C_1 \otimes C_2.$
\end{itemize}
\end{proposition}

We finish this subsection with a note on the co-degree of
the convolution of two rational probability density
functions.
Note that for two independent random variables, $X_1$, $X_2$ say,
it holds that
\[
\Ex |X_1+X_2|^r < \infty \mbox{ if and only if }
\Ex |X_1|^r < \infty  \mbox { and }
\Ex |X_2|^r < \infty.
\]
see e.g. \cite{Rohatgi76}, Problem~4.6.11. Together with
\ref{prop:Ex}  this implies that the co-degree of the convolution
of two rational densities is equal to the minimum of the
co-degrees of these two densities. This fact is used in the text
to track the co-degrees of the conditional probability density
functions arising in the filter.

\newpage
\begin{table}[ht]
\begin{tabular}{ccrr}
&
&
\multicolumn{1}{c}{$\sigma =0.5$}&
\multicolumn{1}{c}{$\sigma =1$}\\
\hline
\( \Ex |Y_{t}| \)&
\( a=0.5 \)& \texttt{0.7202}& \texttt{0.7809}\\
&
\( a=0.9 \)& \texttt{0.7797}& \texttt{1.0279}\\
\hline
\( \mathrm{Var}(|Y_{t}|) \)&
\( a=0.5 \)& \texttt{0.8506}& \texttt{1.3303}\\
&
\( a=0.9 \)& \texttt{1.2994}& \texttt{4.3120}\\
\hline
\( \mathrm{Corr}(|Y_{t+1}|,|Y_{t}|) \)&
\( a=0.5 \)& \texttt{0.0209}& \texttt{0.0619}\\
&
\( a=0.9 \)& \texttt{0.1133}& \texttt{0.2270}\\
\end{tabular}
\caption{Moments of the process $|Y_t|$.
\label{tab:moments}}
\end{table}

~\blfootnote{This table shows the moments of the absoulute
values of the outputs $|Y_t|$ for some
parameter values $a,\Psi=1,\sigma$.}

\newpage
\begin{table}[ht]
\texttt{\begin{tabular}{cccrr}
&
&
&
\multicolumn{1}{c}{$\sigma =0.5$}&
\multicolumn{1}{c}{$\sigma =1$}\\
\hline
\( \hat{a}-a \)&
mean&
\( a=0.5 \)& \texttt{-0.3154}& \texttt{0.0297}\\
&
&
\( a=0.9 \)& \texttt{-0.0522}& \texttt{-0.0000}\\
\cline{2-2} \cline{3-3} \cline{4-4} \cline{5-5}
\multicolumn{1}{c}{}&
std&
\( a=0.5 \)& \texttt{0.2203}& \texttt{0.2321}\\
&
&
\( a=0.9 \)& \texttt{0.2322}& \texttt{0.0591}\\
\hline \( \hat{\Psi }-\Psi  \)& mean&
\( a=0.5 \)& \texttt{-0.0347}& \texttt{0.0074}\\
&
&
\( a=0.9 \)& \texttt{-0.0350}& \texttt{-0.0644}\\
\cline{2-2} \cline{3-3} \cline{4-4} \cline{5-5}
\multicolumn{1}{c}{}&
std&
\( a=0.5 \)& \texttt{0.0642}& \texttt{0.0726}\\
&
&
\( a=0.9 \)& \texttt{0.1010}& \texttt{0.4151}\\
\hline
\( \hat{\sigma }-\sigma  \)&
mean&
\( a=0.5 \)& \texttt{-0.1343}& \texttt{-0.2232}\\
&
&
\( a=0.9 \)& \texttt{-0.0428}& \texttt{0.0142}\\
\cline{2-2} \cline{3-3} \cline{4-4} \cline{5-5}
\multicolumn{1}{c}{}&
std&
\( a=0.5 \)& \texttt{0.4956}& \texttt{0.4751}\\
&
&
\( a=0.9 \)& \texttt{0.3886}& \texttt{0.4262}\\
\end{tabular}}
\caption{Simulation results for the MM estimator.
\label{tab:simresults}}
\end{table}
~\blfootnote{This table shows the results for \( 1000 \)
simulation runs for (simulated) data series of length \( T=1000 \).
Both the mean estimation error (mean) and the standard deviation (std) over
these 1000 simulation are shown in dependence of the true parameters.}

\newpage
Figure \ref{fig:yen_dollar2}:\\
\includegraphics[width=0.8\textwidth]{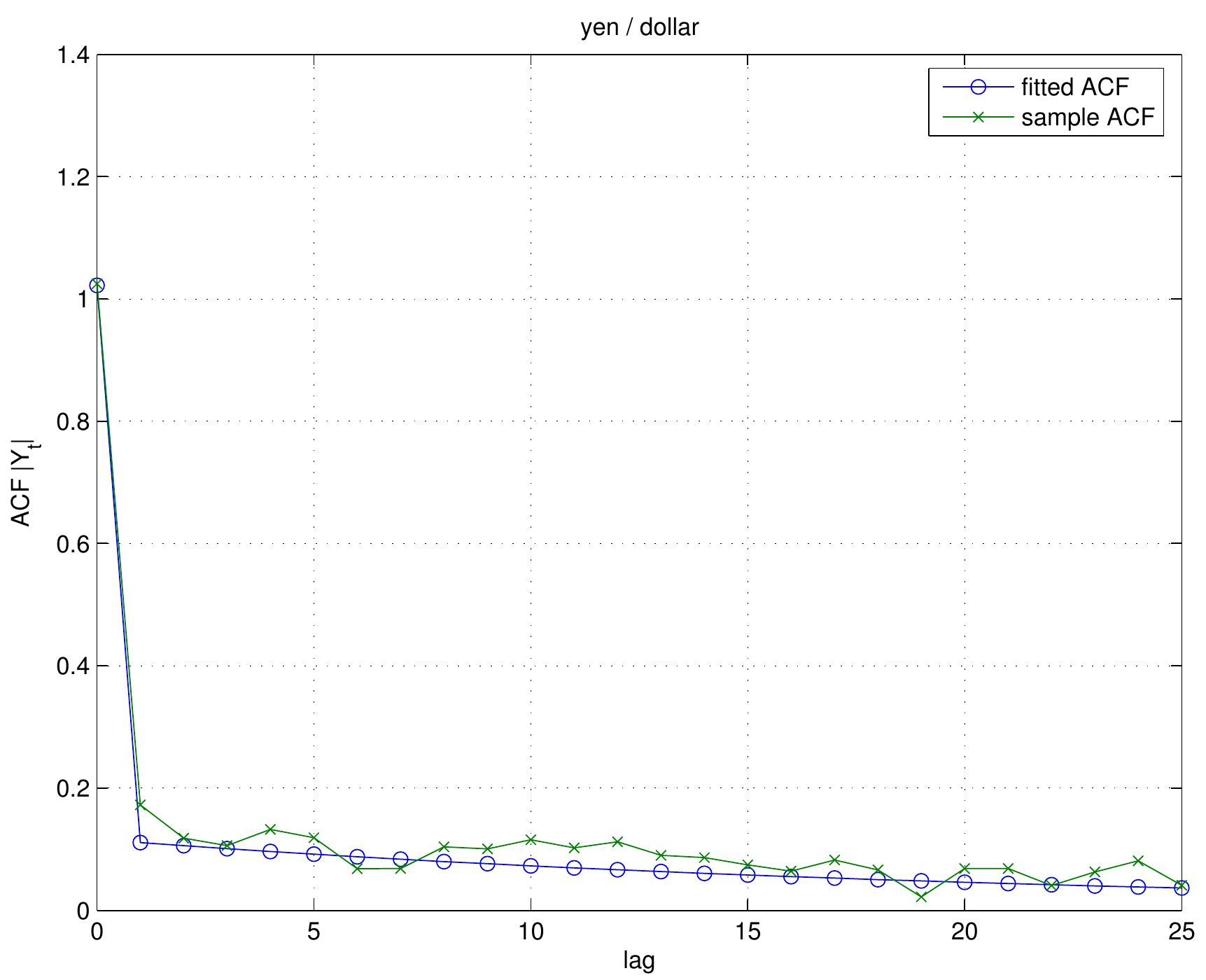}
\begin{figure}[ht]
\caption{Dollar/Yen exchange rate: sample ACF (green) and fitted
ACF of the absolute values.\label{fig:yen_dollar2}}
\end{figure}

\newpage
Figure \ref{fig:yen_dollar3}:\\
\includegraphics[width=0.8\textwidth]{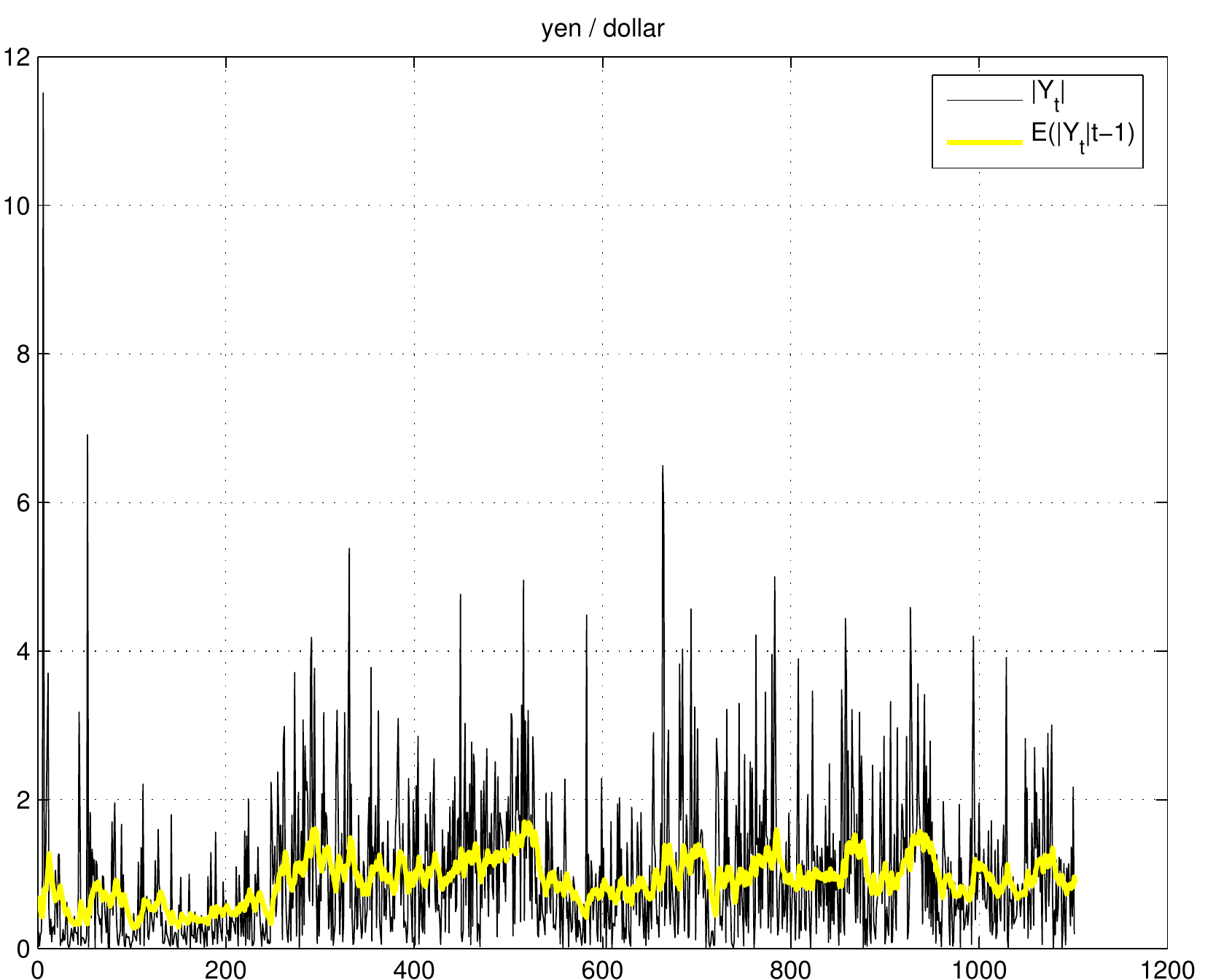}
\begin{figure}[ht]
\caption{Dollar/Yen exchange rate: absolute values of the exchange
rates (black) and the corresponding one step ahead predictions as
given by the filter (yellow).\label{fig:yen_dollar3}}
\end{figure}

\newpage
Figure \ref{fig:simulatedData_fig3}:\\
\includegraphics[width=0.8\textwidth]{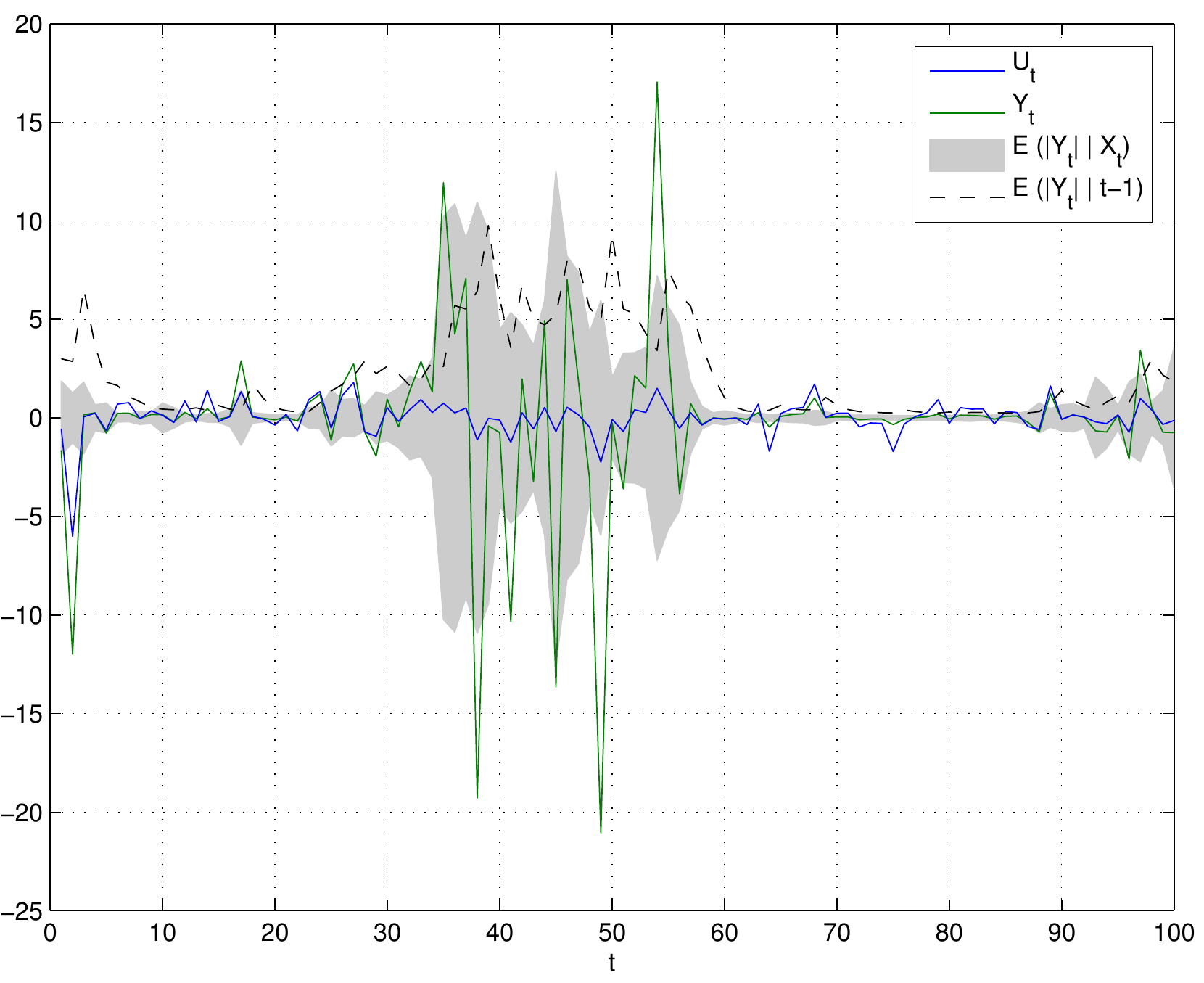}
\begin{figure}[ht]
\caption{Simulated data: time series plot of the simulated
trajectories of the noise process $U_t$ (blue) and of the outputs
$Y_t=V(X_t)U_t$ (green). The gray shaded area is bounded by $\pm
V(X_t)\Ex |U_t|$ i.e. by the conditional expectation of the
absolute values of $Y_t$ given the state $X_t$. The dashed black
line shows the corresponding conditional expectation of $|Y_t|$
given the past observations as computed by the filter, (i.e. the
one step ahead forecasts of
$|Y_t|$).\label{fig:simulatedData_fig3}}
\end{figure}

\newpage
Figure \ref{fig:simulatedData_fig6}:\\
\includegraphics[width=0.8\textwidth]{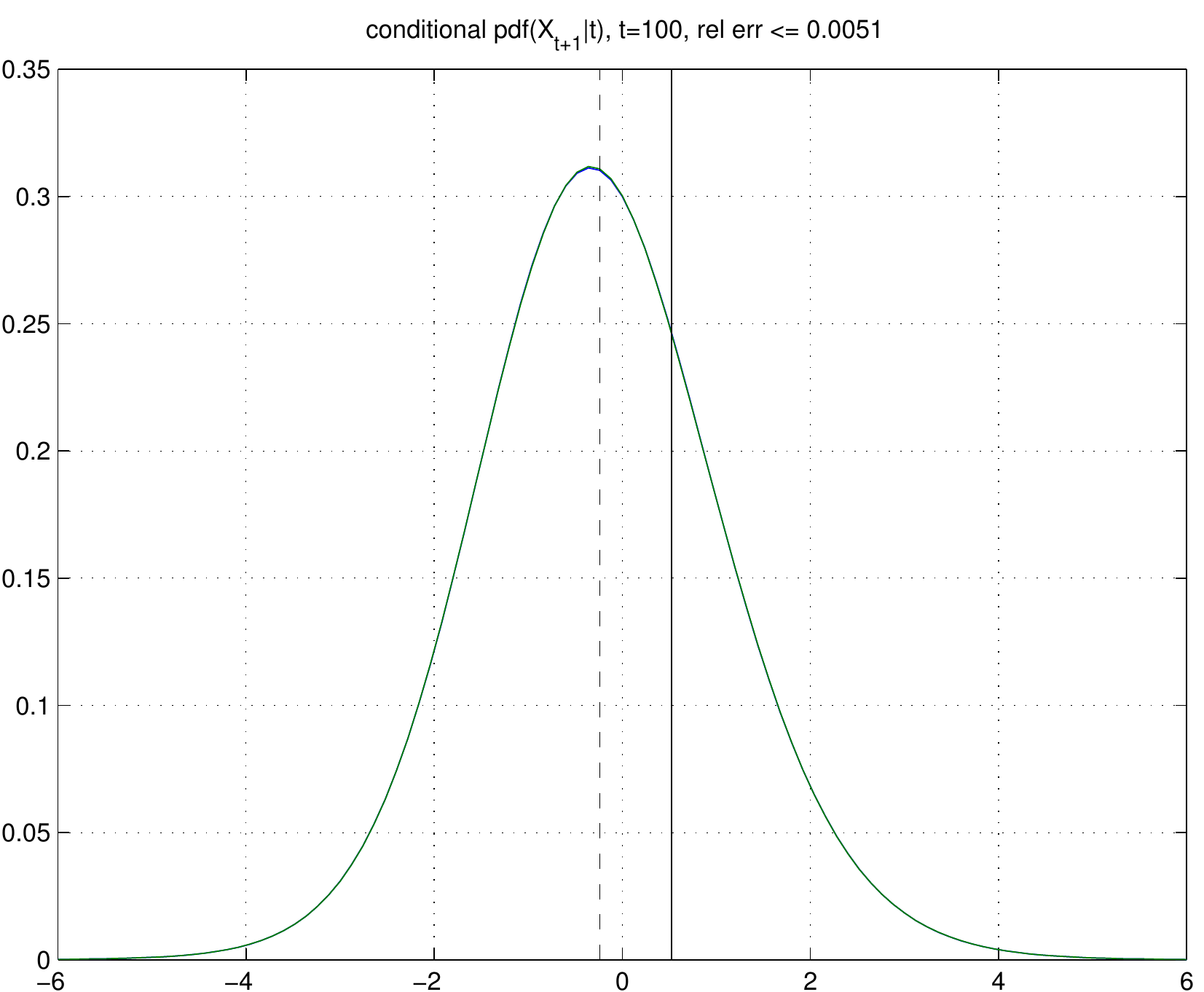}
\begin{figure}[ht]
\caption{Simulated data: This figure shows the conditional pdf
$p_{X_t|Y_1^t}$ for $t=100$. Note that also the approximant
probability density function, $\hat{p}_{X_{t+1}|Y_1^{t}}$ say, as
computed by the positive real balanced truncation method is
plotted. However on this scale the full order probability density
function and the low order approximation can  hardly be
distinguished, since for the relative approximation error
$|p_{X_t|Y_1^t}(x)-\hat{p}_{X_{t+1}|Y_1^{t}}(x)|/p_{X_t|Y_1^t}(x)
\leq 0.0051$ holds by~(\ref{eq:greenbound3}). The state space
dimension of the realization of the spectral summand of
$p_{X_{t+1}|Y_1^{t}}$ is $n=85$ and the spectral summand of
$\hat{p}_{X_{t+1}|Y_1^{t}}$ has order $m=9$. The co-degree of the
corresponding spectral densities
$\Phi_{X_{t+1}|Y_1^{t}}(ix)=p_{X_{t+1}|Y_1^{t}}(x)$ and
$\hat{\Phi}_{X_{t+1}|Y_1^{t}}(ix)=\hat{p}_{X_{t+1}|Y_1^{t}}(x)$ is
$10$.\newline The vertical black line marks the true value
$X_{t+1}$ and the dashed black line marks the corresponding
estimate, i.e. the conditional expectation $\Ex
(X_{t+1}|Y_1^{t})$. \label{fig:simulatedData_fig6}}
\end{figure}

\newpage
Figure \ref{fig:simulatedData_fig7}:\\
\includegraphics[width=0.8\textwidth]{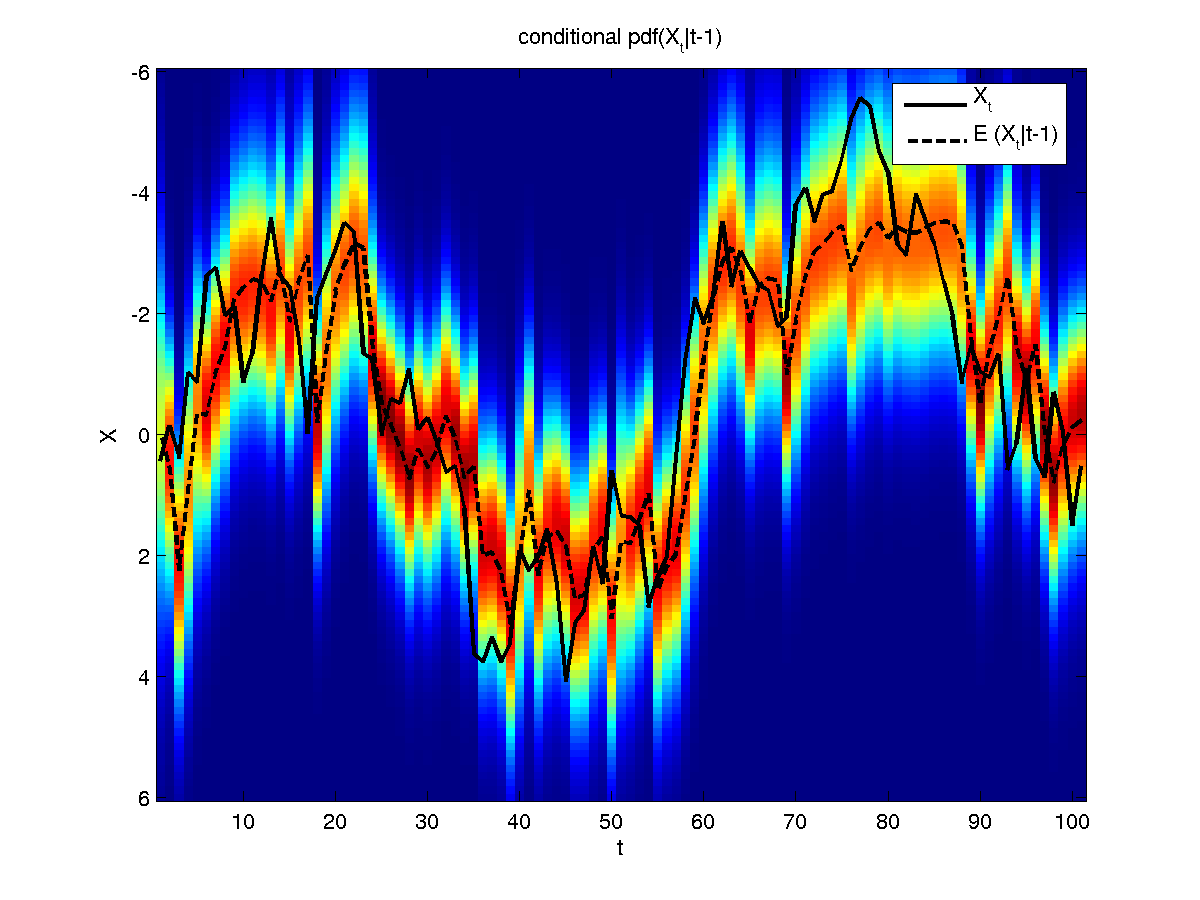}
\begin{figure}[ht]
\caption{Simulated data: This plot shows the evolution of the
conditional densities $p_{X_{t+1}|Y_1^{t}}$. Each ``column'' shows
the conditional $p_{X_{t+1}|Y_1^{t}}$ for a given time $t$, where
high values are coded with red and low values of this pdf are
coded with blue. In addition the solid black line shows the
trajectory of $X_t$ and the blue line marks the corresponding one
step ahead predictions, i.e. the mean of the conditional densities
$p_{X_{t+1}|Y_1^{t}}$.\label{fig:simulatedData_fig7}}
\end{figure}


\begin{thebibliography}{99}
\bibitem{BrigoHanzon98} D. Brigo and B. Hanzon,
  {\em On some filtering problems arising in mathematical
          finance}, Insurance: Mathematics and Economics, vol. 22, 1998,
  pp. 53-64.
\bibitem{CampbellLoMacKinlay97} J. Y. Campbell, A. W. Lo and
 A. C. MacKinlay, {\em The Econometrics of Financial Markets},
 Princeton University Press, Princeton New Jersey,
 2nd ed, 1997.
\bibitem{DesaiPal84} U. B. Desai and D. Pal,
  {\em A transformation approach to stochastic model reduction},
  IEEE Transactions on Automatic Control, vol. AC-29, nr. 12,
  December 1984,pp. 1097-1100.
\bibitem{Engle82} R.F. Engle,
  {\em Autoregressive Conditional Heteroscedasticity with
          Estimates of the Variance of the U.K. Inflation},
          Econometrica,
  vol. 50, 1982, pp. 987-1008.
\bibitem{Faurre76} P. L. Faurre, {\em Stochastic Realization
  Algorithms},pp. 1--25 in: R.K. Mehra and D.G. Lainiotis, {\em
  System Identification: Advances and Case Studies},
  Academic Press, New York, 1976.
\bibitem{GolubVanLoan89} G. Golub and C. VanLoan,
  {\em Matrix Computations,}
  John Hopkins University Press, Maryland, 2nd ed., 1989.
\bibitem{Gourieroux97} C. Gourieroux,
  {\em ARCH Models and Financial Applications}, Springer, New York,
  1997.
\bibitem{Green88AC} M. Green,
  {\em A relative error bound for balanced stochastic truncation},
  IEEE Transactions of Automatic Control, vol. AC-33, nr. 10, pp.
  961-965, October 1988.
\bibitem{Green88LAA} M. Green, {\em Balanced Stochastic Realizations},
  Linear Algebra and its Applications,
  vol. 98, 1988, pp. 211-247.
\bibitem{HanzonOber01} B. Hanzon and R.J. Ober,
  {\em A State-Space Calculus for Rational Probability
  Density Functions and Applications to
  Non-Gaussian Filtering}, SIAM J. Control and
  Optimization, vol. 40, nr.3, 2001, pp. 724-740.
\bibitem{HarveyRuizShephard94} A. Harvey, E. Ruiz and N. Shephard,
  {\em Multivariate Stochastic Variance Models}, Review of Economic
  Studies, vol. 61, 1994, pp. 247-264.
\bibitem{LancasterTismenetsky85} P. Lancaster and M. Tismenetsky,
  {\em The Theory of Matrices}, Academic Press, Orlando, Florida,
  1985.
\bibitem{Lucas96} A. Lucas, {\em Outlier robust unit root analysis},
  Thesis Publishers, Amsterdam, 1996.
\bibitem{MahieuSchotman98} R. Mahieu and P. Schotman,
 {\em An Empirical Application of Stochastic Volatility
 Models}, Journal of Applied Econometrics, vol. 13, June
 1998, pp. 333-360.
\bibitem{Mandelbrot63} B. Mandelbrot,
  {\em The variation of certain speculative prices}, J. Business,
 vol. 36, 1963, pp. 394-419.
\bibitem{Rohatgi76}
  V.K. Rohatgi, {\em An Introduction to Probability Theory and
 Mathematical Statistics}, John Wiley \& Sons, New York, 1976.
 \bibitem{Neuts101} M.F. Neuts, {\em Matrix-Geometric Solutions in
 Stochastic Models: An Algorithmic Approach,} The Johns Hopkins
 University Press, Baltimore, 1981.
 \bibitem{FackrellThesis} M.W. Fackrell, {\em Characterization of
 Matrix-exponential Distributions,} PhD thesis, School of Applied
 Mathematics, Adelaide, 2003.\newline
 http://thesis.library.adelaide.edu.au/uploads/approved/adt-SUA20051207.123257/public/02whole.pdf
\bibitem{Rugh96} Wilson J. Rugh, {\em Linear system theory},
 Prentice-Hall, Upper Saddle River, NJ, 2nd ed., 1996.
\bibitem{Taylor86} S. Taylor, {\em Modelling Financial Time Series},
 John Wiley and Sons, London, 1986.
\bibitem{VanDooren81} P.M. Van Dooren,
 {\em The Generalized Eigenstructure Problem in Linear System
 Theory},
 IEEE Transactions on Automatic Control,
  vol. AC-26, nr. 1, February 1981, pp. 111-129.
\end{thebibliography}
\end{document}